\pdfoutput=1
\RequirePackage{ifpdf}
\ifpdf % We are running pdfTeX in pdf mode
\documentclass[pdftex]{sigma}
\else
\documentclass{sigma}
\fi

\usepackage{bbm}

\newcommand{\One}{\mathbbmss{1}}

\newcommand {\Cee} {{\mathbb C}}

\newcommand {\Zee} {{\mathbb Z}}

\newcommand {\fc} {{\mathfrak{c}}}

\newcommand {\fder} {{\mathfrak{der}}} %

\newcommand {\fg} {{\mathfrak{g}}} %
\newcommand {\fgl} {{\mathfrak{gl}}} %
\newcommand {\fh} {{\mathfrak{h}}}

\newcommand {\fii} {{\mathfrak{i}}} %
\newcommand {\fk} {{\mathfrak{k}}}

\newcommand {\fpe} {{\mathfrak{pe}}} %

\newcommand {\fpgl} {{\mathfrak{pgl}}}

\newcommand {\fpsl} {{\mathfrak{psl}}}

\newcommand {\fq} {{\mathfrak{q}}}

\newcommand {\fs} {{\mathfrak{s}}}

\newcommand {\fsl} {{\mathfrak{sl}}}

\newcommand {\fsvect} {{\mathfrak{svect}}}

\newcommand {\fv} {{\mathfrak{v}}} %
\newcommand {\fvect} {{\mathfrak{vect}}} %

\let\cal\relax
\newcommand {\cal} {\mathcal}

\newcommand {\cL} {{\cal L}}

\newcommand {\cO} {{\cal O}}

\def \opname#1#2%
 {\expandafter\newcommand \csname #1\endcsname {\mathop{\mathrm{#2}}\nolimits }}

\newcommand{\rmname}[1]
 {\expandafter\newcommand \csname #1\endcsname {\mathop{\mathrm{#1}}\nolimits}}

\newcommand{\rmnameii}[2]
 {\expandafter\newcommand \csname #1\endcsname {\mathop{\text{\rm #2}}\nolimits}}

\rmname{act}
\rmname{Ad}
\rmname{Add}
\rmname{ad}
\rmname{Alt}
\rmname{alt}
\rmname{Ann}
\rmname{antidiag}
\rmname{Ber}
\rmname{ber}
\rmname{Bil}
\rmname{Br}
\rmname{But}
\rmname{card}
\rmname{ch}
\rmnameii{Char}{char}
\rmname{Ch}
\rmname{cem}
\rmname{cj}
\rmname{Cliff}
\rmname{cntr}
\rmname{Coker}
\rmname{codim}
\rmname{coind}
\rmname{Conn}
\rmname{conj}
\rmname{const}
\rmname{col}
\rmname{cork}
\rmname{Covect}
\rmname{cpr}
\rmname{diag}
\rmnameii{Div}{div}
\rmname{Def}%Is used by smart.sty; definition in smart blocked now (VM)
\rmname{Der}
\rmname{Diff}
\rmname{Dim}
\rmname{End}
\rmname{Ev}
\rmname{Even}
\rmname{Ext}
\rmname{Gal}
\rmname{gr}
\rmname{grad}
\rmname{Herm}
\rmname{Hom}
\rmname{HT}
\rmnameii{Ht}{ht}
\rmname{hwt}
\rmname{Id}
\rmname{id}
\rmname{ind}
\rmname{Ind}
\rmname{inv}
\rmname{Inv}
\rmname{Coind}%VM
\rmname{Inf}
\rmname{irr}
\rmname{Ker}
\rmname{Le}
\rmname{Lie}
\rmname{lwt}
\rmname{Mat}
\rmname{mult}
\rmname{Mor}
\rmname{nm}
\rmname{Ob}
\rmname{Odd}
\rmname{ord}
\rmname{Osc}
\rmname{per}
\rmname{pet}
\rmname{Pf}
\rmname{Pic}
\rmname{pr}
\rmname{pro}
\rmname{Prime}
\rmname{Proj}
\rmname{prt}
\rmname{pt}
\rmname{ptr}
\rmname{Q}
\rmname{QE}
\rmname{QSt}
\rmname{qch}
\rmname{qet}
\rmname{qtr}
\rmname{rd}
\rmname{rk}
\rmname{row}
\rmname{Reg}
\rmname{Res}
\rmname{salt}
\rmname{Sch}
\rmname{sch}
\rmname{SBr}
\rmname{scalar}
\rmname{sdim}
\rmname{Ser}
\rmname{sign}
\rmname{SIGN}
\rmname{Smbl}
\rmname{spin}
\rmname{srk}
\rmname{ssym}
\rmname{sym}
\rmname{str}

\rmname{st}
\rmname{sgn}
\rmname{sq}
\rmname{symm}
\rmname{supp}
\rmname{Supp}
\rmname{St}
\rmname{Spec}
\rmname{Spm}
\rmname{tr}
\rmname{Vect}
\rmname{vpt}
\rmname{weyl}
\rmname{Weyl}
\rmname{Witt}

\opname{vvol} {{v\hspace{-0.1ex}o\hspace{-0.02ex}l\/}}
\opname{pnt} {\text{\normalfont pt}}
\opname{Span} {{Span}}
\opname{slim} {\overline{\lim}}
\opname{Vol} {{V\hspace{-0.55ex}o\hspace{-0.02ex}l\/}}
\opname{QVol} {{Q\hspace{-0.3ex}V\hspace{-0.55ex}o\hspace{-0.02ex}l\/}}
\opname{PoVol}{{P\hspace{-0.35ex}o\hspace{-0.25ex}V\hspace{-0.55ex}o\hspace{-0.02ex}l\/}}
\opname{BVol} {{B\hspace{-0.2ex}V\hspace{-0.55ex}o\hspace{-0.02ex}l\/}}
\opname{Par} {{P\hspace{-0.3ex}a\hspace{-0.05ex}r\/}}
\opname{Pty} {{P\hspace{-0.3ex}t\hspace{-0.05ex}y\/}}

\rmnameii {IM} {Im}
\rmnameii {RE} {Re}

\opname{Aut} {{A\hspace{-0.2ex}u\hspace{-0.1ex}t\/}}
\opname{GL} {{G\hspace{-0.3ex}L}}
\opname{SL} {{S\hspace{-0.3ex}L}}
\opname{Exp} {{E\hspace{-0.2ex}x\hspace{-0.1ex}p\/}}
\opname{GQ} {{G\hspace{-0.2ex}Q}}
\opname{OSp} {{O\hspace{-0.25ex}S\hspace{-0.15ex}p\/}}
\opname{Out} {{O\hspace{-0.25ex}u\hspace{-0.15ex}t\/}}
\opname{Spp} {{S\hspace{-0.2ex}p\/}}
\opname{SpO} {{S\hspace{-0.2ex}p\hspace{-0.02ex}O\/}}
\opname{Pe} {{P\hspace{-0.25ex}e\/}}
\opname{PGL} {{P\hspace{-0.25ex}G\hspace{-0.25ex}L\/}}
\opname{SPe} {{S\hspace{-0.25ex}P\hspace{-0.25ex}e\/}}
\opname{Spin} {{S\hspace{-0.25ex}p\hspace{-0.05ex}i\hspace{-0.1ex}n\/}}
\opname{Iso} {{I\hspace{-0.25ex}s\hspace{-0.1ex}o\/}}
\opname{Int} {{I\hspace{-0.25ex}n\hspace{-0.1ex}t\/}}
\opname{SSPe} {{S\hspace{-0.25ex}S\hspace{-0.15ex}P\hspace{-0.25ex}e\/}}
\opname{PeU} {{P\hspace{-0.25ex}e\hspace{-0.1ex}U\/}}
\opname{QU} {{Q\hspace{-0.15ex}U\/}}
\opname{U} {{U\/}}

\opname{cGQ} {{\cal G \hspace{-0.2em} Q \/}}
\opname{cSL} {{\cal S \hspace{-0.2em} L \/}}
\opname{cGL} {{\cal G \hspace{-0.2em} L \/}}\opname{cGr} {{\cal G \hspace{-0.2em} r \/}}
\opname{cOSp} {{\cal O \hspace{-0.2em} S \hspace{-0.3em} \it p\/}}
\opname{cPe} {{\cal P \hspace{-1.5pt} \it e\/}}
\opname{cVect} {{\cal V \hspace{-1.5pt} \it e\hspace{-0.1ex}c\hspace{-0.1ex}t\/}}
\opname{cVol} {{\cal V \hspace{-1.5pt} \it o\hspace{-0.1ex}l\/}}
\opname{cAut} {{\cal A \hspace{-0.2em} \it u\hspace{-0.1em}t\/}}
\opname{cCovect} {{\cal C \hspace{-1.5pt}
 \it o\hspace{-0.1ex}v\hspace{-0.1ex}e\hspace{-0.1ex}c\hspace{-0.1ex}t\/}}
\opname{CW} {{C\hspace{-0.15ex}W}}

\newcommand {\ev} {{\bar0}}
\newcommand {\od} {{\bar1}}

\newcommand {\tto} {\longrightarrow}

\opname {Ab} {{\sf Ab}}
\opname {Algs} {{\sf Algs}}
\opname {ASch} {{\sf Aff\;Sch}}
\opname {SSch} {{\sf SSch}}
\opname {Coh} {{\sf Coh}}\opname {Funct} {{\sf Funct}}
\opname {Gr} {{\sf Gr}}
\opname {Grf} {{{\sf Gr}_f}}
\opname {Man} {{\sf Man}}
\opname {MMan} {{\sf MMan}}
\opname {Mods} {{\sf Mods}}
\opname {SMods} {{\sf SMods}}
\opname {Rings} {{\sf Rings}}
\opname {Salgs} {{\sf Salgs}}
\opname {Sets} {{\sf Sets}}
\opname {SSMan} {{\sf SMan}}\opname {Sh} {{\sf Sh}}
\opname {Top} {{\sf Top}}
\opname {Vebun} {{\sf Vebun}}
\opname {Bun} {{\sf Bun}}

\newcommand {\fby}{{\mathfrak{by}}}

\newcommand{\un}{\underline{N}}

\newcommand {\fkl}{{\widetilde{\mathfrak{qv}}}}

\newcommand{\del}{\partial}

\newcommand{\ur}{{\underline{r}}}

\newcommand {\eb}{\boldsymbol{e}}

\newcommand{\veps}{\varepsilon}

\newcommand{\bo}[2]{o^{\{#1,#2\}}}

\numberwithin{equation}{section}

\newtheorem{Theorem}{Theorem}[section]
\newtheorem{Corollary}[Theorem]{Corollary}
\newtheorem{Lemma}[Theorem]{Lemma}

\newtheorem{Problem}[Theorem]{Problem}
\newtheorem{Conjecture}[Theorem]{Conjecture}
\newtheorem{Statement}[Theorem]{Statement}
 { \theoremstyle{definition}

\newtheorem*{Convention}{Convention}
\newtheorem{Example}[Theorem]{Example}
\newtheorem{Remark}[Theorem]{Remark} }

\begin{document}

\allowdisplaybreaks

\newcommand{\arXivNumber}{1711.00638}

\renewcommand{\PaperNumber}{130}

\FirstPageHeading

\ShortArticleName{On Gradings Modulo 2 of Simple Lie Algebras in Characteristic 2}

\ArticleName{On Gradings Modulo 2 of Simple Lie Algebras\\ in Characteristic 2}

\Author{Andrey KRUTOV~$^{\dag\ddag}$ and Alexei LEBEDEV~$^\S$}

\AuthorNameForHeading{A.~Krutov and A.~Lebedev}

\Address{$^\dag$~Institute of Mathematics, Polish Academy of Sciences,\\
\hphantom{$^\dag$}~ul.~\'{S}niadeckich 8, 00-656 Warszawa, Poland}
\Address{$^\ddag$~Independent University of Moscow, Bolshoi Vlasyevskij Pereulok 11, 119002, Moscow, Russia}
\EmailD{\href{mailto:a.o.krutov@gmail.com}{a.o.krutov@gmail.com}}

\Address{$^\S$~Equa Simulation AB, Stockholm, Sweden}
\EmailD{\href{mailto:alexeylalexeyl@mail.Ru}{alexeylalexeyl@mail.ru}}

\ArticleDates{Received January 10, 2018, in final form November 30, 2018; Published online December 10, 2018}

\Abstract{The ground field in the text is of characteristic~2. The classification of modulo~2 gradings of simple Lie algebras is vital for the classification of simple finite-dimensional Lie superalgebras: with each grading, a simple Lie superalgebra is associated, see arXiv:1407.1695. No classification of gradings was known for any type of simple Lie algebras, bar restricted Jacobson--Witt algebras (i.e., the first derived of the Lie algebras of vector fields with truncated polynomials as coefficients) on not less than 3 indeterminates. Here we completely describe gradings modulo 2 for several series of Lie algebras and their simple relatives: of special linear series, its projectivizations, and projectivizations of the derived Lie algebras of two inequivalent orthogonal series (except for ${\mathfrak{o}}_\Pi(8)$). The classification of gradings is new, but all of the corresponding superizations are known. For the simple derived Zassenhaus algebras of height $n>1$, there is an $(n-2)$-para\-metric family of modulo 2 gradings; all but one of the corresponding simple Lie superalgebras are new. Our classification also proves non-triviality of a deformation of a simple $3|2$-dimensional Lie superalgebra (new result).}


\Classification{17B50; 17B20; 17B70}

\section{Introduction}\label{Sintro}

\subsection{Basic definitions}\label{SBackground} Hereafter ${\mathbb K}$ is an algebraically closed field of characteristic $p=2$, unless otherwise specified; all algebras are finite-dimensional; for a review of simple vectorial Lie (super)algebras over ${\mathbb K}$ and basic background, see~\cite{BGLLS2}.

{\bf Lie superlagebras in characteristic 2.} A \textit{Lie superalgebra} is a superspace $\fg=\fg_\ev\oplus\fg_\od$ such that the even part~$\fg_\ev$ is a~Lie algebra, the odd part~$\fg_\od$ is a $\fg_\ev$-module, and on~$\fg_\od$, a squaring (roughly speaking, the halved bracket) is defined as a map
\begin{gather*}
 x\mapsto x^2\quad \text{such that $(ax)^2 = a^2 x^2$ for any $x\in\fg_\od$ and $a\in{\mathbb K}$, and}\\
 \text{$(x+y)^2 - x^2 - y^2$ is a bilinear form on~$\fg_\od$ with values in~$\fg_\ev$.}
\end{gather*}
Then the bracket of odd elements is defined to be
\begin{gather*}
 [x,y] := (x+y)^2 - x^2 - y^2.
\end{gather*}
The Jacobi identity involving odd elements takes the following form:
\begin{gather*}
 \big[x^2,y\big]=[x,[x,y]] \qquad \text{for any} \ x\in\fg_\od, \ y\in\fg_\ev,\qquad
 \big[x^2,x\big]=0\qquad \text{for any} \ x\in\fg_\od.%\label{JI}
\end{gather*}

{\bf Divided powers.} There are two natural integer bases of the commutative algebra $\Cee[x]$ of polynomials in $m$ indeterminates $x=(x_1,\ldots,x_m)$: the monomial one and the basis of divided powers, constructed as follows. For any multi-index $\ur=(r_1,\ldots,r_m)$, where $r_1,\ldots,r_m$ are non-negative integers, we set
\begin{gather*}
u_i^{(r_i)}:=\frac{x_i^{r_i}}{r_i!}\qquad\text{and}\qquad u^{(\ur)}:=\prod_{1\leq i\leq m} u_i^{(r_i)}.
\end{gather*}
These $u^{(\ur)}$ form an integer basis of~$\Cee[x]$. Their multiplication relations are
\begin{gather}\label{mRel}
u^{(\ur)}\cdot u^{(\underline{s})} = \binom{\ur + \underline{s}}{\ur} u^{(\ur + \underline{s})},\qquad \text{where}\quad\binom{\ur + \underline{s}}{\ur} = \prod_{1\leq i\leq m} \binom{r_i + s_i}{r_i}.
\end{gather}

Being interested in simple Lie algebras, observe that every simple $\Zee$-graded Lie algebra~$\fg\!=\!\oplus\fg_i$ of vector fields is \textit{transitive}, i.e., such that
\begin{gather*}
 [\fg_{\leq 0}, x]=0 \ \ \text{for a given $x\in \fg_{>0}$ implies $x=0$, where $\fg_{\leq0}=\oplus_{i\leq0}\fg_i$, $\fg_{>0} = \oplus_{i>0}\fg_i$.}
\end{gather*}

For an arbitrary field~${\mathbb K}$ of characteristic~$p>0$, it is easy to see that the Lie algebra $\fder \, {\mathbb K}[x]$ is not transitive in the $\Zee$-grading induced by the standard $\Zee$-grading of~${\mathbb K}[x]$, i.e., $\deg x_i=1$ for all $i$. The situation is remedied if we consider the commutative algebra~${\mathbb K}[u]$ spanned by all the elements~$u^{(\ur)}$ with multiplication relations~\eqref{mRel}. For any $m$-tuple $\un = (N_1,\ldots,N_m)$, where $N_i$ are either positive integers, or infinity, denote (we set
$p^{\infty}:=\infty$)
\begin{gather*}
 \cO(m;\un) := {\mathbb K}[u;\un] = \Span_{\mathbb K}\big( u^{(\ur)} \,|\, r_i < p^{N_i} \big).
\end{gather*}
The algebra~${\mathbb K}[u]$ and its subalgebras~${\mathbb K}[u;\un]$ are called the \textit{algebras of divided powers}; they are analogs of the polynomial algebra. Let $\One:=(1,\dots,1)$ denote the shearing vector $\un$ with the smallest values of heights $N_i$ of the indeterminates.

Clearly, if $N_i\neq 1$ for at least one $i$, then $\cO(m;\un)$ has more than $m$ generators: namely \mbox{$y_{i,j}:=u_i^{(p^{j-1})}$}. Any derivation~$D$ of a given algebra is determined by the values of~$D$ on the generators, so~$\fder\,\cO(m;\un)$ has more than $m$ analogs of partial derivations: one for each generator, with a~functional parameter $f_{ij}(u)\in \cO(m;\un)$~-- coefficient of each partial derivative~$\partial_{y_{i,j}}$, whereas the coefficients of $\partial_{y_{i,j}}\in\fvect\big(\sum N_i;\One\big)$ belong to $\cO\big(\sum N_i;\One\big)$. Following the definitions over~$\Cee$ we are interested in subalgebras of~$\fder\,\cO(m;\un)$ such that the dimension of the nonpositive part (in the standard $\Zee$-grading) is equal to~$m$.

Solution was found long ago: one has to introduce \textit{distinguished} partial derivatives $\partial_i$, each of them serving as several partial derivatives (corresponding to $u_i$, $u_i^{(p)}$, $u_i^{(p^2)}$, \dots) at once (if instead of the usual powers we use divided ones, the sign $(-1)^{j-1}$ is not needed):
\begin{gather*}
\partial_i\big(u_j^{(k)}\big):=\delta_{ij}u_j^{(k-1)}\ \ \text{ for all $k$, i.e., }\partial_i=\sum_{j\geq 1} (-1)^{j-1}y_{i,1}^{p-1}\cdots y_{i,j-1}^{p-1}\partial_{y_{i,j}}.
\end{gather*}
The \textit{general vectorial Lie algebra of distinguished derivations}\footnote{For $p>0$, the Lie algebra $\fvect(m;\un)$ is called \textit{Jacobson--Witt} algebra if $\un\neq \One$; it is usually denoted $W(m;\un)$ for any $m$ and $\un$; if $m=1$, it is called \textit{Zassenhaus} algebra. Jacobson--Witt algebras are simple if $m>1$; there is no special name for the simple derived $\fvect^{(1)}(1;\underline{n})$, see a discussion in~\cite{GZ}. The Lie algebra of \textit{divergence-free}, or ``special'', vector fields is denoted $\fsvect(m;\un)$, usually abbreviated to $S(m;\un)$.}{\samepage
\begin{gather*}
 \fvect(m;\un) = \Span_{\mathbb K} ( f \del_k \,|\, f \in \cO(m;\un), \, k=1,\dots, m)
\end{gather*}
is simple (except for $m=1$ and $p=2$).}

{\bf $\boldsymbol{p}$-structure.} If $\fg$ is a Lie algebra, then for every $x\in\fg$, the operator $(\ad_x)^p$ is a derivation of~$\fg$. If this derivation is inner for every $x\in\fg$, then the Lie algebra $\fg$ is said to be \textit{restricted} or having a $p$-\textit{structure}. More specifically, a $p$-\textit{structure} on $\fg$ is a map $[p]\colon \fg\to\fg$, $x\mapsto x^{[p]}$ such that
\begin{gather}
 \big[x^{[p]}, y\big] = (\ad_x)^p(y)\qquad\text{for any $x,y\in\fg$},\\
(ax)^{[p]}=a^px^{[p]}\qquad \text{for any~}a\in{\mathbb K},~x\in\fg,\nonumber\\
(x+y)^{[p]}=x^{[p]}+y^{[p]}+ \sum \limits_{1\leq i\leq p-1}s_i(x, y) \qquad\text{for any~}x,y\in\fg,\label{restricted-3}
\end{gather}
where $s_i(x, y)$ is the coefficient of $\lambda^{i-1}$ in $(\ad_{\lambda x+y})^{p-1}(x)$.

\begin{Remark}\quad
\begin{enumerate}\itemsep=0pt
\item[1)] If the Lie algebra $\fg$ is without center, then the last two conditions of~\eqref{restricted-3} follow from the first one. There might be more than one $p$-structure on one Lie algebra; all of them are equal modulo center. Hence, on any simple Lie algebra, there is at most one $p$-structure.
\item[2)] The following condition is sufficient for a Lie algebra $\fg$ to possess a $p$-structure: for a basis $\{g_i\}_{i\in I}$ of $\fg$, there exist elements $g_i^{[p]}$ such that
\begin{gather*} \big[g_i^{[p]}, y\big]=(\ad_{g_i})^{p}(y)\qquad \text{for any} \quad y\in\fg.
\end{gather*}
\end{enumerate}
\end{Remark}

We consider the following problem.
\begin{Problem}\label{Prob}
For any finitely generated commutative group $G$, classify $G$-gradings of simple finite-dimensional Lie algebras over~${\mathbb K}$.
\end{Problem}

\textbf{Lie algebras for which a solution of Problem~\ref{Prob} is known}.
Although we are interested in a particular case of Problem~\ref{Prob}, let us briefly review the known general results.
\begin{itemize}\itemsep=0pt
\item {\it For $p=0$}, for a clear exposition of the solution (for the $\Zee$- and $\Zee/n$-gradings), see the
book~\cite[around p.~500]{H}, and also~\cite{EK}.

\item {\it For $p\neq 2$}, \cite{Ko} is a very lucid review of the cases where the solution of Problem~\ref{Prob} has been found; for further details, see~\cite{BK,KPS}; for examples of applications of certain $G$-gradings, see~\cite{Kos}.

\item {\it For $p=2$}, the result of~\cite{BK} is as follows: all $G$-gradings of the Jacobson--Witt algebra~$W(m;\One )$ (it is $\fvect(m;\One )$ in our notation) for $m\geq3$ are given by $G$-gradings of the corresponding algebra of divided powers~$\mathcal{O}(m;\One )$ due to an isomorphism of their automorphism group schemes~\cite{Sk01}. The classification of such gradings is given. Any $G$-grading of $\cO(m;\One )$ is equivalent, up to an algebra automorphism, to one that can be described as follows. For a given~$s$ such that $0\leq s \leq m$, and $a_1,\ldots,a_m\in G$, set
\begin{gather}\label{BK}
\cO_g = \Span\bigg\{ (1+x_1)^{j_1}\cdots(1+x_s)^{j_s}x_{s+1}^{j_{s+1}}\cdots x_{m}^{j_m} \,|\, j_i= 0 \ \text{or} \ 1, \ \sum_{1\leq i\leq m} j_i a_i= g \bigg\},
\end{gather}
where $a_1,\ldots,a_m\in G$ are the respective degrees of the indeterminates generating $\mathcal{O}(m;\One )$, i.e.,
\begin{gather*}
1+x_1,\ldots, {1+x_s}, x_{s+1},\ldots, x_m.
\end{gather*}
\end{itemize}

For $p>0$ and the \textit{restricted} Lie algebras considered in \cite{BK} (for $p=2$: Jacobson--Witt algebras~$\fvect(m;\un)$ for $m\geq3$; for $p>2$: Jacobson--Witt algebras, simple relatives of Hamiltonian algebras and algebras of divergence-free vector fields), the only possible $\Zee/2$-gradings are generated by those described in equation~\eqref{BK}. For \textit{non-restricted} Lie algebras, we know several examples of $\Zee/2$-gradings of Kaplansky algebras, see~\cite{BLLS}.

\subsection[Problem~\ref{Prob} for $G=\Zee/2$ and $p=2$: an application of the solution]{Problem~\ref{Prob} for $\boldsymbol{G=\Zee/2}$ and $\boldsymbol{p=2}$: an application of the solution}\label{Sz/2grad}

In \cite{BLLS1}, there are offered two methods for constructing a~simple finite-dimensional Lie superalgebra from every simple finite-dimensional Lie algebra \textit{over a field of characteristic~$2$}; it is proved that
every simple finite-dimensional Lie superalgebra can be obtained by one of these two methods (queerification and ``method~2''). The ``method 2'' depends on the $\Zee/2$-gradings of the simple Lie algebra which is being superized. Let us recall the method:

Let $\fg=\fg_{\ev}\oplus\fg_{\od}$ be a simple Lie algebra with a $\Zee/2$-grading $\gr$. Let the 1-\textit{step restricted closure} of $\fg$ associated with the grading $\gr$ be
\begin{gather}
\fg^{\langle 1\rangle}:=\text{the minimal Lie subalgebra of the restricted closure $\overline{\fg}$ containing $\fg$}\nonumber\\
\hphantom{\fg^{\langle 1\rangle}:=}{} \ \text{and all the elements $x^{[2]}$, where $x\in\fg_{\od}$.} \label{resCl}
\end{gather}
Clearly, there is a single way to extend the grading $\gr$ from $\fg$ to $\fg^{\langle 1\rangle}$; we assume this extension performed. On the space of $\fg^{\langle 1\rangle}$, define the structure of a Lie superalgebra denoted (in what follows, we often omit indicating the grading $\gr$ since it enters the definition of $\fg^{\langle 1\rangle}$)
\begin{gather}\label{S(g)}
\fs\big(\fg^{\langle 1\rangle}, \gr\big)\qquad \text{by setting $x^2:=x^{[2]}$ for any $x\in\fg_{\od}$,}
\end{gather}
and retaining the bracket of any even element with any other element. The Lie superalgebra~$\fs\big(\fg^{\langle 1\rangle}\big)$ is simple, see~\cite{BLLS1}. Our strategic goal is classification of simple Lie (super)algebras, so we formulate our description of $\Zee/2$-gradings $\gr$ of $\fg$ in terms of superizations $\fs\big(\fg^{\langle 1\rangle}\big)$ whenever we can.

\subsection[How we seek $\Zee/2$-gradings]{How we seek $\boldsymbol{\Zee/2}$-gradings}
Let $\fg$ be a Lie algebra, and let $\fg = \fg_{\ev}\oplus \fg_{\od}$ be its $\Zee/2$-grading. For an any $x\in\fg$, we denote by~$x_{\ev}$ its even part and by $x_{\od}$ its odd part. We have
\begin{gather*}
[x,y] = [x_{\ev} + x_{\od}, y_{\ev} + y_{\od}] = [x_{\ev}, y_{\ev}] + [x_{\od}, y_{\od}] + [x_{\ev}, y_{\od}] + [x_{\od},y_{\ev}] \qquad \text{for any $x,y\in\fg$}.
\end{gather*}
This implies
\begin{gather*}
([x,y])_{\ev} = [x_{\ev}, y_{\ev}] + [x_{\od}, y_{\od}]\qquad\text{and}\qquad ([x,y])_{\od} = [x_{\ev}, y_{\od}] + [x_{\od}, y_{\ev}].
\end{gather*}
Finally, let $U \in\End(\fg)$, where $\fg$ is considered as a vector space, be the projection to the~$\fg_\od$, i.e., $Ux = x_{\od}$ and $(I-U)x = x_{\ev}$, where $I$ is the identity operator. A~linear operator $U$ is such a projection to the odd part of $\fg$ in some $\Zee/2$-grading if and only if it satisfies the following conditions for all $x,y\in\fg$:
\begin{subequations}\label{cond}
 \begin{gather}
 U^2= U,\label{cond1nonlin}\\
 U[x,y] = [Ux,(I - U)y] + [(I - U)x,Uy] = -2[Ux,Uy] + [Ux,y] + [x,Uy].\label{cond1lin}
\end{gather}
\end{subequations} Having fixed a basis in $\fg$, we express these conditions in terms of structure constants $c^{ij}_k$ as (summation over repeated indices is assumed)
\begin{subequations}\label{cond2}
\begin{gather}
U^k_jU^i_k= U^i_j,\label{cond2nonlin}\\
c^{ij}_lU^l_k= c^{il}_kU^j_l+ c^{lj}_kU^i_l- 2c^{lm}_kU^i_lU^j_m.\label{cond2lin}
\end{gather}
\end{subequations}
For $p=2$, equation~\eqref{cond2lin} becomes linear and easy (for example, for \textit{Mathematica}-based computer package \textit{SuperLie}, see~\cite{Gr}) to solve. There remains, however, a~problem to be solved:{\samepage
%\begin{subequations}\label{cond3}
\begin{gather*}
\text{condition \eqref{cond2nonlin} is still quadratic},\\ %\label{cond3nonlin}\\
\text{we need equivalence classes of solutions $U$ mod $\Aut(\fg)$, not individual operators}. %\label{cond3lin}
 \end{gather*}}
%\end{subequations}

{\bf A comment: gradings and derivations.} %\label{grader}
For a general discussion of relation between gradings and derivations, and interesting examples for $p=2$, see~\cite{BLLS}. For $p=3$, an example illustrating the said discussion is the $\Zee/4$-grading of the Skryabin
algebra~$\fby$, see~\cite{GL}.

The Lie algebra $\fder\, \fvect(m;\un)$ of all derivations of $\fvect(m;\un)$ coincides with the $p$-envelope of $\fvect(m;\un)$. So the problem of classification of equivalence classes of $\Zee/p$-gradings of $\fvect(m;\un)$ reduces to the classification of equivalence classes of toral elements in $\fder\, \fvect(m;\un)$. The equivalence classes of maximal tori with respect to the group of automorphisms are described in~\cite{Kuz, T}, and the answer depends on $N_1-2$ parameters for $m=1$.

\subsection{Our results}
\begin{enumerate}\itemsep=0pt
\item[1)] We classify the $\Zee/2$-gradings of the Lie algebras of series $\fsl$, their simple subquotients, and simple derived or subquotients of both orthogonal series~-- ${\mathfrak{o}}_I$ and ${\mathfrak{o}}_\Pi$,~-- except for~${\mathfrak{o}}_\Pi(8)$. In all the cases considered, these gradings yield the known superizations of the corresponding Lie algebras.
Note that one of the gradings of ${\mathfrak{o}}^{(1)}(3)$ has no analogs among the $\Zee/2$-gradings of the simple 3-dimensional Lie algebras for $p\neq 2$. This unusual grading has analogs among $\Zee/2$-gradings
of~$\fvect^{(1)}(1;\underline{n})$, of which ${\mathfrak{o}}^{(1)}(3) \simeq \fvect^{(1)}(1;\underline{2})$ is a~particular case.

\item[2)] The classification of $\Zee/2$-gradings of ${\mathfrak{o}}^{(1)}(3)$ has one more application: it gives the shor\-test known proof of the fact that the deform of ${\mathfrak{oo}}_{I\Pi}^{(1)}(1|2)$, one of the two superizations of~${\mathfrak{o}}^{(1)}(3)$, found in \cite{BGL1}, is a ``true'' one, not ``semitrivial'', see Remark~\ref{R_oII}; this fact is new and unexpected (although the deform itself, ${\mathfrak{oo}}_{II}^{(1)}(1|2)$, is known).

\item[3)] We describe an $(n-2)$-parametric collection of $\Zee/2$-gradings of $\fvect^{(1)}(1;\underline{n})$. For $n=2$, these gradings yield three non-isomorphic Lie superalgebras.

For $n>2$, these $\Zee/2$-grading yield
\begin{enumerate}\itemsep=0pt
\item[a)] purely even Lie superalgebra~$\fvect^{(1)}(1;\underline{n}|0)$;
\item[b)] an $(n-2)$-parametric family of filtered deforms of $\fk(1;\underline{n-1}|1)$;
\item[c)] a filtered deform of~$\fq(\fvect(1;\underline{n-1}))$.
\end{enumerate}
\end{enumerate}

{\bf Examples where $\boldsymbol{\Zee/2}$-gradings of ``the same algebra'' for $\boldsymbol{p\neq 2}$ and $\boldsymbol{p=2}$ differ or where previously unknown $\boldsymbol{\Zee/2}$-gradings have been found.} These are the most interesting of our results.
\begin{enumerate}\itemsep=0pt
\item[1)] It is known that for any $p$, there is only one simple Lie algebra of dimension 3: for $p\neq 2$, it is ${\mathfrak{o}}(3)\simeq\fsl(2)$; for $p=2$, it is ${\mathfrak{o}}^{(1)}(3)$, the derived of ${\mathfrak{o}}(3)$.

For $p\neq 2$, it is also known that there is only one (as always, up to an automorphism) nontrivial $\Zee/2$-grading, and hence if the superization procedure had been defined for $p\neq 2$, the superization of $\fsl(2)$ would have been unique, $\fsl(1|1)$.

For $p=2$, there are 2 inequivalent nontrivial $\Zee/2$-gradings of ${\mathfrak{o}}^{(1)}(3)$, as we will see below. The corresponding non-isomorphic, see~\cite{Leb}, Lie superalgebras are ${\mathfrak{oo}}^{(1)}_{I\Pi}(1|2)$ and ${\mathfrak{oo}}^{(1)}_{II}(1|2)$. (It is an interesting open problem to find out if the Lie superalgebra ${\mathfrak{o}}{\mathfrak{o}}_{II}^{(1)}(1|2n)$ is a deform of~${\mathfrak{o}}{\mathfrak{o}}_{I\Pi}^{(1)}(1|2n)$ for $n>1$; this requires to consider $\Zee/2$-gradings and deforms of~${\fk(2n+1;\One)}$ for~$n>0$.)

\item[2)] For $p\neq 2$, the $\Zee/2$-gradings of the Zassenhaus algebra $\fvect(1;\underline{n})$ were only known for $n=1$, i.e., for the restricted algebras, see \cite{BK} (and, for $n>1$, obvious $\Zee/2$-gradings with some of either $x$ or $x+1$ declared odd). We have found new $(n-2)$-parametric family of {$\Zee/2$-gra}\-dings in the non-restricted case. These gradings yield superizations which have no known analogs in the case $p\neq 2$.
\end{enumerate}

\begin{Remark}\label{Perm} Describing $\fder\, \fg$ for various Lie algebras $\fg$, we heavily rely on Permyakov's results~\cite{Per} used in~\cite{BGLL2}. Consider, for example, $\fsl(2n)$ for $n>1$. Actually, Permyakov \cite{Per} does not explicitly describe $\fder\, \fsl(2n)$, he just shows that ${\dim H^1(\fsl(2n), \fsl(2n))=1}$, this is the codimension of the subalgebra of inner derivations in $\fder\, \fsl(2n)$. However, it is easy to see that every element of $\fpgl(2n)$ described here represents a derivation of $\fsl(2n)$; different elements represent different derivations; every inner derivation is represented by some element, and the codimension of the subalgebra or inner derivations of $\fsl(2n)$ is equal to $1$; therefore, $\fpgl(2n)$ is isomorphic to $\fder\, \fsl(2n)$. Similar considerations apply when we refer to~\cite{BGLL2, Per} while describing algebras of derivations of other Lie algebras.
\end{Remark}

\subsection{Open questions}
For $G=\Zee/2$, the gradings \eqref{BK} with $s=0$ of $\fvect(m;\One)$ for $m>2$ yield $\fvect(k;\One|m-k)$, where the vectors $\One$ have $m$ and $k$ coordinates, respectively. These are known superizations. If $s\neq 0$ in \eqref{BK} and $a_i= \od$ for some of the indices $i\leq s$, then the corresponding indeterminates $1+x_i$ are odd; at the moment we are unable to identify the resulting Lie superalgebra.

The list of simple Lie algebras for which Problem~\ref{Prob} is open:
\begin{itemize}\itemsep=0pt
\item For $p>3$, the simple vectorial Lie algebras for the shearing vector $\un\neq\One $, and the deforms (results of deformations) thereof (for their classification and description, see \cite{Kos,SkH, S}).
\item For $p=3$, the simple vectorial Lie algebras; in particular, exceptional ones, mainly discovered by Skryabin and lucidly described in~\cite{GL}, and their deforms to be described (for a review, see~\cite{BGL1}); the deforms of ${\mathfrak{o}}(5)$ and Brown algebras described in~\cite{BLW}.
\item For $p=2$, there are many examples of simple Lie algebras of types not existing for $p\neq 2$, see~\cite{BGLL2,BGLLS2, BGLLS1,BLLS,Ei,GZ, SkT1}, and \cite{BGL1, BGL}.
\item We also have to consider the simple relative of ${\mathfrak{o}}_\Pi(8)$.
\end{itemize}

\begin{Conjecture}[D.~Leites]\label{2conj}
%\sssec{A conjecture (D.~Leites)} \label{2conj}
For $($the simple derived of$)$ Lie algebras with indecomposable Cartan matrix not considered in this note, i.e., parametric families of Weisfeiler--Kac algebras, all their $\Zee/2$-gradings correspond, conjecturally, to their known superizations listed in~{\rm \cite{BGL}}.
\end{Conjecture}

\section[The $\fsl(n)$ and $\fpsl(2n)$ for $n>2$: same answer as for $p\neq 2$]{The $\boldsymbol{\fsl(n)}$ and $\boldsymbol{\fpsl(2n)}$ for $\boldsymbol{n>2}$: same answer as for $\boldsymbol{p\neq 2}$}\label{S_sl}

Clearly, $\fsl(0|n)=\fsl(n)$ and $\fpsl(0|2n)=\fpsl(2n)$.

\begin{Theorem}\label{Tsl} All $\Zee/2$-gradings of the Lie algebras $\fsl(n)$ for $n>2$ and $\fpsl(2n)$ for $n>2$ are analogous to those for $p=0$, i.e., correspond to $\fsl(k|n-k)$ for $k=0,\ldots,\big[\frac{n}{2}\big]$ and $\fpsl(k|2n-k)$ for $k=0,\ldots,n$.
\end{Theorem}

\begin{Remark} The algebra~$\fpsl(2)$ is not simple, so we do not consider algebras~$\fsl(2)$ and~$\fpsl(2)$. We do not consider the algebra~$\fpsl(4)$ in this section either, since its algebra of derivations is more complicated than the algebra of derivations of~$\fpsl(2n)$ for $n>2$; however, $\fpsl(4)$ is isomorphic to~${\mathfrak{o}}^{(2)}_\Pi(6)$, and the answer in this case is given in Remark~\ref{psl4}.
\end{Remark}

\begin{proof}
{\it The $\fsl(2n+1)$ series $(n>0)$.} All derivations of $\fsl(2n+1)$ are inner ones, see~\cite{Per}, so the general solution of the linear equation \eqref{cond1lin} is $U=\ad_A$, where $A\in \fsl(2n+1)$. The nonlinear equation~\eqref{cond1nonlin} reduces then to $\ad_A = (\ad_A)^2$. Now recall that there is the 2-structure on $\fsl(2n+1)$ given by $A^{[2]} =A^2$, see \cite{BGL,BLLS1}. Therefore
\begin{gather*}
\ad_A = (\ad_A)^2= \ad_{A^{[2]}}= \ad_{A^2},
\end{gather*}
and since $\fsl(2n+1)$ has no center, we have $A = A^2$, i.e., $A$ is a projection. Since $\tr A$ is equal to the dimension of the subspace onto which~$A$ projects (just look at the normal form of the matrix of the projection), it follows that
\begin{gather}
\text{a given projection $A$ belongs to $\fsl$ if and only if it is a projection}\nonumber\\ \text{onto a subspace of even dimension}.\label{tr0}
\end{gather}

If $A\in\fsl(2n+1)$ is a projection, then the Lie superalgebra obtained from it by ``method 2'' of~\cite{BLLS1} is isomorphic to
\begin{gather*}
\fsl(\dim \Ker A \, |\dim \IM A)\cong \fsl(\dim \IM A\, |\dim \Ker A).
\end{gather*}
So we see there are $n+1$ equivalence classes of $\Zee/2$-gradings of $\fsl(2n+1)$ and $n+1$ of its nonisomorphic superizations (including the trivial purely even one). We can enumerate them either as
\begin{gather*}
\fsl(2n+1-2k|2k)\cong \fsl(2k|2n+1-2k), \qquad \text{where $k =0,\dots,n$,}
\end{gather*}
or, more simply, as $\fsl(k|2n+1-k)$, where $k = 0,\dots,n$. So the answer is the same as for $p\neq 2$.

{\it The $\fsl(2n)$ series $(n\geq 2)$.} The Lie algebra $\fder\ \fsl(2n)$ can be identified, see \cite{Per}, with $\fpgl(2n)$ in the sense that for any $D\in\fder\ \fsl(2n)$, there is ${A_D\in\fgl(2n)}$ such that $D$
coincides with the restriction of $\ad_{A_D}$ to $\fsl(2n)$.

These elements $A_D$ are defined modulo center; in particular, one can take $A_{D^2}$ to be $(A_D)^2$. Therefore, $D$ satisfies the condition \eqref{cond1nonlin} if and only if $(A_D)^2 = A_D + cI_{2n}$ for some $c\in{\mathbb K}$. Let $d\in{\mathbb K}$ be a root of the equation $d^2=d+c$; set $A'_D = A_D + dI_{2n}$, then $(A'_D)^2 = A'_D$. So we see that an operator $U$ on $\fsl(2n)$ satisfies the conditions \eqref{cond1nonlin} if and only if it can be represented as a restriction of $\ad_A$ to $\fsl(2n)$ for some projection $A\in\fgl(2n)$, and then the Lie superalgebra obtained from $U$ by ``method 2'' is isomorphic to
\begin{gather*}
{\fsl(\dim \Ker A\, |\dim \IM A)\cong \fsl(\dim \IM A\, |\dim \Ker A)}.
\end{gather*}
So there are $n+1$ nonisomorphic superizations of $\fsl(2n)$ (including the trivial purely even one): $\fsl(k|2n-k)$, where $k = 0,\dots,n$. Again, the answer is the same as for $p\neq 2$.

{\it The $\fpsl(2n)$ series, $n>2$.} The algebra $\fder\ \fpsl(2n)$ is isomorphic to $\fder\ \fsl(2n)$ when $n>2$, see~\cite{Per}. So the arguments from the previous subsection apply. There are $n+1$ non-isomorphic superizations of $\fpsl(2n)$ (including the purely even one), which can be enumerated as \smash{$\fpsl(k|2n-k)$}, where $k = 0,\dots,n$.
\end{proof}

\begin{Remark}\label{WG} A $\Zee/2$-grading of $\fgl(2n+1)$ for $n>2$ can produce a Lie superalgebra which is not isomorphic to any superalgebra of the form $\fgl(k|2n+1-k)$. More specifically, such a superalgebra would be isomorphic to the direct sum of a Lie superalgebra of the form ${\fsl(k|2n+1-k)}$, where $0\leq k\leq n$, and $0|1$-dimensional center (i.e., 1-dimensional odd center). Namely, represent $\fgl(2n+1)$ as $\fsl(2n+1)\oplus\fc$; since it is a direct sum, you can introduce gradings on the summands independently. On $\fsl(2n+1)$, introduce some grading which produces ${\fsl(k|2n+1-k)}$, then declare the center odd. The resulting algebra is not isomorphic to $\fgl(k|2n+1-k)$ for any $k$, since the center of $\fgl(k|2n+1-k)$ is even. But since this grading does not result in any new \textit{simple} Lie superalgebra, we do not consider gradings of $\fgl$ here.
\end{Remark}

\section[The ${\mathfrak{o}}^{(1)}(2n+1)$ series]{The $\boldsymbol{{\mathfrak{o}}^{(1)}(2n+1)}$ series}\label{So2n+1}

Clearly, ${\mathfrak{o}}{\mathfrak{o}}_{B_1B_2}(0|n)={\mathfrak{o}}_{B_2}(n)$ and ${\mathfrak{o}}{\mathfrak{o}}_{B_1B_2}(n|0)={\mathfrak{o}}_{B_1}(n)$.

\begin{Example}[${\mathfrak{o}}^{(1)}(3)$] Observe the occasional isomorphism ${\mathfrak{o}}^{(1)}(3)\simeq\fvect^{(1)}(1;\underline{2})$, so this case is considered in more details in Section~\ref{ssW2}. In~\cite{BGL}, the deforms of the superization of ${\mathfrak{o}}^{(1)}(3)$ associated with one of the three $\Zee/2$-gradings of ${\mathfrak{o}}^{(1)}(3)$ are described; but until now it was unclear if these deforms are \textit{true} ones, i.e., the deformed algebra is not isomorphic to the initial algebra, as is the case for \textit{semitrivial} deformations corresponding to certain integrable cocycles from nontrivial cohomology class, for examples, see \cite{BLLS, Ri}. The results of Section~\ref{ssW2} prove that one of the deforms found in \cite{BGL} is a true one.
\end{Example}

\begin{Theorem}\label{T_oOdd} For any $n\geq 1$, all $\Zee/2$-gradings of the Lie algebra ${\mathfrak{o}}^{(1)}(2n+1)$ correspond to Lie superalgebras ${\mathfrak{oo}}^{(1)}_{II}(2n+1-2k|2k)$, where $0\leq k\leq n$ $($or, which is the same but looks simpler, ${\mathfrak{oo}}^{(1)}_{II}(k|2n+1-k)$, where $0\leq k\leq n)$, and ${\mathfrak{oo}}^{(1)}_{I\Pi}(2n+1-2k|2k)$, where $1\leq k\leq n$.
\end{Theorem}

\begin{proof}Since any nondegenerate symmetric bilinear form over a space of dimension $2n+1$ is equivalent to $I_{2n+1}$, the algebra of derivations of ${\mathfrak{o}}^{(1)}(2n+1)$ is isomorphic to ${\mathfrak{o}}(2n+1)/\fc$, the quotient of ${\mathfrak{o}}(2n+1)$ modulo center, or the subalgebra of traceless elements of ${\mathfrak{o}}(2n+1)$, see \cite[Statement~3.8.1a]{BGLL2} and Appendix~\ref{appendixA}.

Again, an element $A\in \fder\, {\mathfrak{o}}^{(1)}(2n+1)$ describes a $\Zee/2$-grading of ${\mathfrak{o}}^{(1)}(2n+1)$ if and only if $A^2=A$, i.e., $A$ is a projection.

By definition, an operator $A\in\fgl(V)$ belongs to ${\mathfrak{o}}_B(V)$, where $B$ is a nondegenerate symmetric bilinear form, if and only if
\begin{gather}\label{(*)}
B(Ax,y) + B(x, Ay) = 0 \qquad \text{for all $x,y\in V$.}
\end{gather}
In what follows, we assume that $A$ is a projection. It is easy to check that \eqref{(*)} is automatically satisfied if $x,y\in \IM A$ or if $x,y\in \Ker A$. If $x\in \IM A$ and $y\in \Ker A$, then \eqref{(*)} is equivalent to $B(x,y) = 0$. So $A\in {\mathfrak{o}}_B(V)$ if and only if $\IM A$ and $\Ker A$ are orthogonal with respect to $B$. For~$A$ to be traceless, $\dim\IM A$ must be even, see condition~\eqref{tr0}.

Let $A$ be a projection belonging to $\fder\, {\mathfrak{o}}^{(1)}_B(2n+1)$ with $\dim \IM A = 2k>0$. Denote restrictions of $B$ to $\IM A$ and $\Ker A$ by $B_{\IM A}$ and $B_{\Ker A}$, respectively.

Since $B$ is nondegenerate, $\IM A \oplus \Ker A = V$ and $\IM A \perp \Ker A$ with respect to $B$, it follows that these restrictions are nondegenerate.

Since $\dim \Ker A = 2n+1-2k$ is odd, $B_{\Ker A}$ is equivalent to $I_{2n+1-2k}$. Since $\dim \IM A = 2k$ is even, it follows that $B_{\IM A}$ may be equivalent to either $I_{2k}$, or $\Pi_{2k}$. In the latter case, the
resulting Lie superalgebra is isomorphic to ${\mathfrak{oo}}^{(1)}_{I\Pi}(2n+1-2k|2k)$; in the former case, it is isomorphic to ${\mathfrak{oo}}^{(1)}_{II}(2k|2n+1-2k)$. Note that the collection ``${\mathfrak{oo}}^{(1)}_{II}(2k|2n+1-2k)$, where $0\leq k\leq n$'' can be described more simply as ``${\mathfrak{oo}}^{(1)}_{II}(k|2n+1-k)$, where $0\leq k\leq n$''.

So we get $2n+1$ nonisomorphic Lie superalgebras from ${\mathfrak{o}}^{(1)}(2n+1)$, including the purely even case. To see that all these cases are really realizable with some projection $A$, the following approach can be used. Recall that in an odd-dimensional space, all nondegenerate symmetric forms are equivalent, so for any such form $B$ and any symmetric invertible matrix, there is a~basis in which the matrix of $B$ is equal to the given one. In particular, for a given $k$, there are bases in which the Gram matrix of $B$ is either $I_{2n+1}$ (call it \textit{the first basis}) or $\diag(I_{2n+1-2k}, \Pi_{2k})$ (call it \textit{the second basis} for this value of~$k$).

Take the operator $A$ whose matrix in the \textit{first} basis is $\diag(0_{2n+1-2k}, I_{2k})$; clearly, $A^2=A$ and $\dim \IM A=2k$. The space $\IM A$, which is spanned by the last $2k$ vectors of the basis, is, clearly, orthogonal, relative the form $B$, to $\Ker A$, which is spanned by the first $2n+1-2k$ vectors of the basis. So $A\in \fder\, {\mathfrak{o}}_B$, and hence determines a grading of ${\mathfrak{o}}^{(1)}_B$. It is easy to see that the Lie superalgebra obtained from ${\mathfrak{o}}^{(1)}_B$ by ``method 2'' of~\cite{BLLS1} is ${\mathfrak{oo}}^{(1)}_{II}(2n+1-2k|2k)$.

Analogously, an operator $A$ with the same matrix $\diag(0_{2n+1-2k}, I_{2k})$ in the \textit{second} basis (it does not matter what is the matrix of $A$ in the first basis) determines a grading that corresponds to the Lie superalgebra ${\mathfrak{oo}}^{(1)}_{I\Pi}(2n+1-2k|2k)$.

It is possible to show that there are $2n+1$ equivalence classes of $\Zee/2$-gradings. It could have happened that there were more gradings than nonisomorphic superizations if two inequivalent gradings would have yielded isomorphic Lie superalgebras. This, however, does not happen; we skip the details.\end{proof}

\begin{Remark}\label{WGo} Similarly to $\fgl(2n+1)$, see Remark~\ref{WG}, it is clear that certain $\Zee/2$-gradings of ${\mathfrak{o}}(2n+1)$ can produce Lie superalgebras not isomorphic to any superalgebra of the form
${\mathfrak{oo}}_{II}(2n+1-2k|2k)$ or ${\mathfrak{oo}}_{I\Pi}(2n+1-2k|2k)$. For example, certain $\Zee/2$-gradings can yield a~direct sum of ${\mathfrak{oo}}^{(1)}_{II}(2n+1-2k|2k)$ and $0|1$-dimensional center. But again, such gradings do not produce any new \textit{simple} Lie superalgebras, so, being actually interested in new simple Lie superalgebras, we do not consider such gradings of ${\mathfrak{o}}(2n+1)$ in full generality in this paper.
\end{Remark}

\section[The ${\mathfrak{o}}^{(1)}_{I}(2n)$ series]{The $\boldsymbol{{\mathfrak{o}}^{(1)}_{I}(2n)}$ series}\label{SoI2n}

\begin{Theorem}\label{T_oIEv} For $n>2$, all $\Zee/2$-gradings of the Lie algebra ${\mathfrak{o}}^{(1)}_{I}(2n)$ correspond to Lie superalgebras ${\mathfrak{oo}}^{(1)}_{II}(m|2n-m)$, where $1\leq m\leq n$, and ${\mathfrak{oo}}^{(1)}_{I\Pi}(2n-2k|2k)$, where $0\leq k\leq n-1$.
\end{Theorem}
\begin{proof} If a bilinear form $B$ on vector space $V$ of dimension $2n$ is equivalent to $I_{2n}$, then the algebra $\fder\, {\mathfrak{o}}^{(1)}_B(2n)$ can be identified, see \cite[Statement~3.8.1a]{BGLL2} and Appendix~\ref{appendixA}, with ${\mathfrak{o}}_B(2n)/\fc$ in the sense that for any $D\in\fder\, {\mathfrak{o}}^{(1)}_B(2n)$, there is $A_D\in{\mathfrak{o}}_B(2n)$ such that $D$ coincides with the restriction of $\ad_{A_D}$ to ${\mathfrak{o}}^{(1)}_B(2n)$. These elements $A_D$ are defined up to a central element; in particular, one can take $A_{D^2}$ to be $(A_D)^2$. By arguments similar to the ones we used in the case of $\fsl(2n)$, one
can show that an operator $U$ on ${\mathfrak{o}}^{(1)}_B(2n)$ satisfies the conditions \eqref{cond1nonlin} if and only if it can be represented as a restriction of $\ad_A$ to ${\mathfrak{o}}^{(1)}_B(2n)$ for some projection
$A\in{\mathfrak{o}}_B(2n)$.

By arguments we used in the previous section, a projection $A$ belongs to ${\mathfrak{o}}_B(2n)$ if and only if $\IM A$ and $\Ker A$ are orthogonal with respect to $B$. We denote restrictions of $B$ to $\IM A$ and $\Ker A$ by $B_{\IM A}$ and $B_{\Ker A}$, respectively. These restrictions have to be nondegenerate.

If $\dim \IM A = 2k+1$, where $0\leq k\leq n-1$, then these restrictions are equivalent to $I_{2k+1}$ and $I_{2n-2k-1}$, respectively, and the resulting Lie superalgebra is isomorphic to
\begin{gather*}
{\mathfrak{o}}^{(1)}_{II}(2k+1|2n-2k-1)\simeq {\mathfrak{o}}^{(1)}_{II}(2n-2k-1|2k+1).
\end{gather*}

If $\dim \IM A = 2k$, where $1\leq k\leq n-1$, then $B_{\IM A}$ can be equivalent to either $I_{2k}$, or $\Pi_{2k}$, while $B_{\Ker A}$ can be equivalent to either $I_{2n-2k}$, or $\Pi_{2n-2k}$. However, it is impossible for $B_{\IM A}$ to be equivalent to $\Pi_{2k}$ while $B_{\Ker A}$ is equivalent to $\Pi_{2n-2k}$ at the same time, because the direct sum of these two forms is equivalent to~$\Pi_{2n}$. All the other combinations are possible, and, depending on them, the resulting Lie superalgebra can be isomorphic to either ${\mathfrak{o}}^{(1)}_{II}(2k|2n-2k)$, or ${\mathfrak{o}}^{(1)}_{I\Pi}(2k|2n-2k)$, or ${\mathfrak{o}}^{(1)}_{\Pi I}(2k|2n-2k)\simeq{\mathfrak{o}}^{(1)}_{I\Pi}(2n-2k|2k)$.

If $\dim \IM A = 0$ or $2n$, then the resulting Lie superalgebra is purely even. So, as described above, ${\mathfrak{o}}^{(1)}_B(2n)$ has $2n+1$ nonisomorphic superizations.\end{proof}

\begin{Remark}\label{R_vAut} The algebras ${\mathfrak{o}}^{(1)}_I(2)$ and ${\mathfrak{o}}^{(1)}_I(4)$ are not simple, so we do not consider them here.
\end{Remark}

\section[The simple relatives of ${\mathfrak{o}}_{\Pi}(2n)$ series]{The simple relatives of $\boldsymbol{{\mathfrak{o}}_{\Pi}(2n)}$ series}\label{SoPi}

The algebras ${\mathfrak{o}}_{\Pi}(2)$ and ${\mathfrak{o}}_{\Pi}(4)$ do not have simple relatives, so we do not consider them. We have been unable to classify $\Zee/2$-gradings of the simple relative of ${\mathfrak{o}}_{\Pi}(8)$ so far.

\begin{Convention} We will denote the simple relative of ${\mathfrak{o}}_{\Pi}(2n)$ by ${\mathfrak{o}}^{(2)}_\Pi(2n)/\fc$ for both even and odd values of $n$, keeping in mind that the center might be trivial. Similarly, we will denote the simple relatives of ${\mathfrak{oo}}_{\Pi\Pi}(2k|2n-2k)$ and $\fpe^{(2)}(n)$ by ${\mathfrak{oo}}^{(2)}_{\Pi\Pi}(2k|2n-2k)/\fc$ and $\fpe^{(2)}(n)/\fc$, respectively, for both even and odd values of $n$.
\end{Convention}

\begin{Theorem}\label{T_oPi} For $n=3$ or $n\geq 5$, all $\Zee/2$-gradings of the simple relative of ${\mathfrak{o}}_{\Pi}(2n)$ $($that is, of ${\mathfrak{o}}^{(2)}_\Pi(2n)$, if $n$ is odd, or of ${\mathfrak{o}}^{(2)}_\Pi(2n)/\fc$, if $n$ is even$)$ correspond to the simple relatives of the corresponding superizations, i.e., ${\mathfrak{oo}}^{(2)}_{\Pi\Pi}(2k|2n-2k)$, where $0\leq k\leq \big\lfloor\frac{n}{2}\big\rfloor$, and $\fpe^{(2)}(n)$, if $n$ is odd, or ${\mathfrak{oo}}^{(2)}_{\Pi\Pi}(2k|2n-2k)/\fc$, where $0\leq k\leq \big\lfloor\frac{n}{2}\big\rfloor$, and $\fpe^{(2)}(n)/\fc$, if $n$ is even.
\end{Theorem}

\begin{proof} If $n$ is odd, then the center of ${\mathfrak{o}}^{(2)}_\Pi(2n)$ is trivial, so in this case ${\mathfrak{o}}^{(2)}_\Pi(2n)/\fc\simeq {\mathfrak{o}}^{(2)}_\Pi(2n)$.

For a symmetric bilinear form $B$ on space $V$, consider the Lie algebra
\begin{gather}
\tilde{\mathfrak{o}}_{B}(V):= \{M\in\fgl(V)\,|\, \text{there is $c_M\in{\mathbb K}$ such that} \nonumber\\
\hphantom{\tilde{\mathfrak{o}}_{B}(V):= \{}{} B(Mx,y) + B(x,My) = c_MB(x,y) \ \text{for all $x,y\in V$}\}.\label{tilde-o}
\end{gather}
As is not difficult to see, for $B\sim\Pi_{2n}$ we have
\begin{gather*}
\tilde{\mathfrak{o}}_{\Pi_{2n}}(V)={\mathfrak{o}}_{\Pi_{2n}}(V) \ltimes {\mathbb K} d_n, \qquad \text{where $d_n = \diag(0_n, 1_n)$ or $\diag(1_n, 0_n)$}.
\end{gather*}
Then, according to \cite{BGLL2, Per}
\begin{gather*}%\label{noDoubt}
\fder\big({\mathfrak{o}}^{(2)}_{\Pi_{2n}}(V)/\fc\big)\simeq\tilde{\mathfrak{o}}_{\Pi_{2n}}(V)/\fc.
\end{gather*}
(To compare with our previous results, observe that if $B\sim I_{2n}$, then $\tilde{\mathfrak{o}}_{I_{2n}}(V) = {\mathfrak{o}}_{I_{2n}}(V)$.)

As before, a given operator $U\in \fder\big({\mathfrak{o}}^{(2)}_B(V)/\fc\big) \simeq \tilde{\mathfrak{o}}_B(V)/\fc$ satisfies the condition \mbox{$U^2=U$} if and only if the corresponding element (equivalence class) of $\tilde{\mathfrak{o}}_B(V)$ contains a~projec\-tion~$A$. Let us consider the values~$c_A$, see~\eqref{tilde-o}, can take.

{$c_A = \ev$}: In this case, the definition \eqref{tilde-o} is equivalent to the statement that $\IM A$ and $\Ker A$ are orthogonal with respect to $B$. It means that the restrictions of $B$ to $\IM A$ and $\Ker A$ have to be nondegenerate. Since $B$ is anti-symmetric (i.e., $B(x,x) = 0$ for any $x\in V$), these restrictions are anti-symmetric as well, which means that dimensions of $\IM A$ and $\Ker A$ are even. The resulting Lie superalgebra in this case is isomorphic to
\begin{gather*}
{\mathfrak{oo}}^{(2)}_{\Pi\Pi}(\dim \IM A\, |\dim\Ker A)/\fc\simeq{\mathfrak{oo}}^{(2)}_{\Pi\Pi}(\dim \Ker
A\, |\dim\IM A)/\fc.
\end{gather*}

{$c_A = \od$}: In this case, the condition \eqref{tilde-o} is equivalent to the statement that both $\IM A$ and $\Ker A$ are isotropic with respect to $B$. Since $\IM A\oplus \Ker A = V$, this means that
\begin{gather*}
\dim \IM A = \dim \Ker A = n
\end{gather*}
and there is an invertible linear map $f\colon \Ker A\to (\IM A)^*$ such that
\begin{gather*}
B(x,y) = (f(y))(x) \qquad \text{for all $x\in\IM A$, $y\in \Ker A$.}
\end{gather*}
The resulting Lie superalgebra in this case is isomorphic to $\fpe^{(2)}(n)/\fc$.

{$c_A\neq \ev,\od$}: This is impossible, because in this case, \eqref{tilde-o} would mean that $\IM A$ and $\Ker A$ are both isotropic and orthogonal to each other with respect to $B$, i.e., $B(x,y) = 0$ for all $x,y\in V$.

So we get $\big\lfloor\frac{n}{2}\big\rfloor+2$ nonisomorphic superizations of ${\mathfrak{o}}^{(2)}_\Pi(2n)/\fc$ (including the purely even one), which are ${\mathfrak{oo}}^{(2)}_{\Pi\Pi}(2k|2n-2k)/\fc$, where $0\leq k\leq
\big\lfloor\frac{n}{2}\big\rfloor$, and $\fpe^{(2)}(n)/\fc$.
\end{proof}

\begin{Remark}\label{psl4} The algebra ${\mathfrak{o}}^{(2)}_\Pi(6)$ is isomorphic to $\fpsl(4)$, $\fh_\Pi^{(1)}(4;\One)$ and $\fsvect^{(1)}(3;\One)$, see~\cite{ChKu}. The fact that it has exactly three superizations implies two facts: 1)~the following superizations of these algebras are isomorphic, and 2)~these algebras have no other non-trivial superizations non-isomorphic to these ones:
\begin{gather*}
 {\mathfrak{oo}}^{(2)}_{\Pi\Pi}(2|4)\simeq\fsvect^{(1)}(0|3)\simeq \fpsl(3|1)\simeq\fsvect^{(1)}(2;\One|1), \qquad \text{here} \ \ \sdim=8|6, \\
\fpe^{(2)}(3)\simeq\fh_\Pi^{(1)}(0|4)\simeq\fpsl(2|2)\simeq\fh_{\Pi\Pi}^{(1)}(2;\One |2)\simeq\fsvect^{(1)}(1;\One|2), \qquad \text{here} \ \ \sdim=6|8.
\end{gather*}
\end{Remark}

\section[The $\fvect^{(1)}(1;\underline{n})$ series]{The $\boldsymbol{\fvect^{(1)}(1;\underline{n})}$ series}\label{S_vect}
We describe the $\Zee/2$-gradings of $\fvect^{(1)}(1;\underline{n-1})$ in Theorem~\ref{T_vect}. First, we need some definitions.

Let $\fv_n$ be the minimal Lie subalgebra of the restricted closure of $\fvect^{(1)}(1;\underline{n-1})$ that contains $\fvect^{(1)}(1;\underline{n-1})$ and all the elements $\mathcal{X}^{[2]}$, where $\mathcal{X}\in\fvect^{(1)}(1;\underline{n-1})$. As a vector space, $\fv_n$ can be written as a direct sum
\begin{gather}\label{fv_n}
 \fv_n = \fvect(1;\underline{n-1}) \oplus {\mathbb K} \del^2.
\end{gather}
Thus, $\dim \fv_n = 2^{n-1}+1$. Observe that $\fv_n = \big(\fvect^{(1)}(1;\underline{n-1})\big)^{\langle 1\rangle}$, see equation~\eqref{resCl}, associated with the $\Zee/2$-grading induced by the standard $\Zee$-grading of $\fvect^{(1)}(1;\underline{n-1})$.

For $n>2$, consider the Lie superalgebra $\fq(\fvect(1;\underline{n-1}))$. As a vector space, it can be written as a direct sum $\fv_n\oplus\Pi(\fvect(1;\underline{n-1}))$, where $\Pi$ is the change of parity
functor. Consider the Lie algebra~$\fvect(1;\underline{n})$ in an indeterminate~$z$ in order to distinguish it from the indeterminate~$x$ in equation~\eqref{evPart}. Let
\begin{gather*}
 X_{-2} = \del^2,\ X_{-1} = \del,\ X_0 = z\del,\ \ldots,\ X_{2^{n-1}-2} = z^{2^{n-1}-1}\partial,
\end{gather*}
be the basis in the even part~$\fv_n$ of~$\fq(\fvect(1;\underline{n-1}))$ and
\begin{gather*}
 Y_{-1} = \Pi \del,\ Y_0 = \Pi (z\del),\ \ldots,\ Y_{2^{n-1}-2} = \Pi \big(z^{(2^{n-1}-1)}\del\big),
\end{gather*}
be the basis in the odd part of~$\fq(\fvect(1;\underline{n\!-\!1}))$. The Lie superalgebra structure in~$\fq(\fvect(1;\underline{n\!-\!1}))$ is given by following formulas
\begin{gather*}
 [X_{-2}, X_k] = \big[\del^2, z^{(k-1)}\del\big] = X_{k-2}\qquad\text{for $k = 1,\ldots,2^{n-1}-2$},\\
 [X_{-2}, Y_k] = \Pi\big[\del^2, z^{(k-1)}\del\big] = Y_{k-2}\qquad\text{for $k=1,\ldots,2^{n-1}-2$},\\
 [X_{-2}, X_{-1}] = [X_{-2}, X_{0}] = [X_{-2}, Y_{-1}] = [X_{-2}, Y_{0}] = [X_{-1}, Y_{-1}] = 0,\\
 [X_k, X_m] = \big[z^{(k-1)}\del, z^{(m-1)}\del\big] = \binom{k+m+2}{k+1}X_{k+m}\\
 \hphantom{[X_k, X_m] =}{} \ \text{for $k,m = - 1,\ldots,2^{n-1}-2$, except for $m=k=-1$,}\\
 [X_k, Y_m] = \Pi\big[z^{(k-1)}\del, z^{(s-1)}\del\big] = \binom{k+m+2}{k+1}Y_{k+m} \\
 \hphantom{[X_k, Y_m] =}{} \ \text{for $k,m = -1,\ldots,2^{n-1}-2$, except for $m=k=-1$,}\\
 [Y_k, Y_m] = \big[z^{(k-1)}\del, z^{(m-1)}\del\big] = \binom{k+m+2}{k+1}X_{k+m}\\
 \hphantom{[Y_k, Y_m] =}{} \ \text{for $k,m=-1,\ldots,2^{n-1}-2$, except for $m=k=-1$,}\\
 (Y_{-1})^2= X_{-2}, \\
 (Y_k)^{2} = \big(z^{(k+1)}\del\big)^{[2]} = \binom{2k + 1}{k} X_{2k} \qquad\text{for $k=0,\ldots,2^{n-1}-2$}.
\end{gather*}

On the superspace of~$\fq(\fvect(1;\underline{n-1}))$, define another Lie superalgebra structure,
call it $\fkl_{n-1}$, by the following formulas:
\begin{gather*}
 [X_{-2}, X_{-1}] := 0,\qquad [X_{-2}, X_0]: = 0,\\
 [X_{-2}, X_m] := X_{m-2}\qquad\text{for~$m = 0,\ldots,2^{n-1}-2$},\\
 [X_k, X_m] := \binom{k+m+2}{k+1} X_{k+m}\qquad\text{for~$m,k = -1,\ldots,2^{n-1}-2$, except for $m=k=-1$,}\\
 [Y_k, Y_m] := \binom{k+m+2}{k+1} X_{k+m}\qquad\text{for~$m,k=-1,\ldots,2^{n-1}-2$, except for $m=k=-1$,}\\
 [X_{-1}, Y_{-1}] := Y_{-1},\qquad [X_{-2},Y_{-1}] = Y_{-1}, \qquad [X_{-2},Y_0] = Y_0,\\
 [X_{-2},Y_m] = Y_m + Y_{m-2}, \qquad m = 1,\ldots,2^{n-1}-2\\
 [X_k, Y_m] := \binom{k+m+2}{k+1} Y_{k+m} + \binom{k+m+2}{k+1} Y_{k+m+1}\\
\hphantom{[X_k, Y_m] :=}{} \ \text{for $k,m =-1,\ldots,2^{n-1}-2$, except for $m=k=-1$,}\\
 (Y_{-1})^2 := X_{-2} + X_0,\qquad (Y_k)^{2} := \binom{2k + 1}{k} X_{2k} \qquad\text{for $k=0,\ldots,2^{n-1}-2$}.
\end{gather*}

I.~Shchepochkina observed that the Lie superalgebra~$\fkl_{n-1}$ is a filtered deform of the Lie superalgebra~$\fq(\fvect(1;\underline{n-1}))$.

\begin{Lemma}Lie superalgebras $\fkl_{n-1}$ and $\fq(\fvect(1;\underline{n-1}))$ are not isomorphic.
\end{Lemma}
\begin{proof} Let us prove that $\fkl_{n-1}$ is not isomorphic to any Lie superalgebra of the form $\fq(\fg)$, where $\fg$ is a Lie algebra. Any superalgebra of the form $\fq(\fg)$ possesses the following property: if $x\in\fq(\fg)_\ev$ acts nilpotently on $\fq(\fg)_\ev$, then it acts nilpotently on $\fq(\fg)_\od$ as well. Indeed, if $(\ad_x)^k|_{\fq(\fg)_\ev} = 0$ for some positive integer $k$, then $(\ad_x)^k \Pi y = \Pi\big( (\ad_x)^k y \big)= 0$ for any $y\in\fg$.

On the other hand, from the relations above one can see that $X_{-1}$ acts nilpotently on $(\fkl_{n-1})_\ev$ (more specifically, $(\ad_{X_{-1}})^{2^{n-1}+1}|_{(\fkl_{n-1})_\ev} = 0$), but $[X_{-1}, Y_{-1}] = Y_{-1}$, so the action of $X_{-1}$ on $(\fkl_{n-1})_\od$ is not nilpotent.
\end{proof}

Consider the Lie algebra~$\fvect(1;\underline{n})$ in an indeterminate~$x$. For a basis of $\fvect(1;\underline{n})$ we take $\eb_i := x^{(i+1)}\del$ for $i=-1,0, \ldots,2^n-2$.

\begin{Theorem}\label{T_vect} All $\Zee/2$-gradings of the Lie algebra $\fvect^{(1)}(1;\underline{n})$ correspond to the one of the following Lie superalgebras:
\begin{enumerate}\itemsep=0pt
\item[$1)$] purely even $\fvect^{(1)}(1;\underline{n}|0)$;
\item[$2)$] $\fk(1; \underline{n-1}|1)$, described in~{\rm \cite{BGLLS2}}, and its $(n-2)$-parametric family of filtered deforms described below in the proof;
\item[$3)$] the Lie superalgebra $\fkl_{n-1}$ for $n>2$.
\end{enumerate}
\end{Theorem}

\begin{proof}Let $u\in \cO(1;\underline{n})$ be a linear combination of only even powers of the indeterminate $x$, and let $a$ be its constant term\footnote{The notation $u(0)$ we use for brevity in what follows is meaningless, strictly speaking, because the divided power polynomials can not be evaluated at any $x\in{\mathbb K}$ for $p>0$; by writing like this we mean $a\in{\mathbb K}\cdot1$ which is the value of $u$ modulo the maximal ideal of $\cO(1;\underline{n})$; this is a standard abuse of notation.}. In other words, let $u,a\in \cO(1;\underline{n})$ be such that
\begin{gather}\label{u-prop}
 x\cdot(\del u) = 0, \qquad \del a = 0, \qquad u^{2} = {}a^2.
\end{gather}

Consider any derivation $D_u \in \fder\ \cO(1;\underline{n})$ of the form
\begin{gather*}
D_u = \bigg(u+x+xu\sum\limits_{1\leq i\leq n-1} a^{2^i - 2}\big(\del^{2^i} u\big)\bigg)\del + \sum\limits_{1\leq i\leq n-1} a^{2^i}\del^{2^i}.
\end{gather*}
Properties \eqref{u-prop} imply (after rather lengthy calculations we omit) that $(D_u)^2 = D_u$, so $D_u$ describes a $\Zee/2$-grading of $\cO(1;\underline{n})$.

Consider the linear map $\ad_{D_u}$ on the Lie algebra $\fder\, \cO(1;\underline{n})$ given by $D\mapsto [D_u,D]$. Due to the Jacobi identity, we see that $\ad_{D_u}\in\fder(\fder\, \cO(1;\underline{n}))$, where the outer $\fder$ is for derivations of the Lie algebra $\fder\, \cO(1;\underline{n})$. Since the Lie algebra $\fder\, \cO(1;\underline{n})$ has a $2$-structure given by $D^{[2]} = D^2$, we have $(\ad_{D_u})^2 = \ad_{D_u}$. Thus, $\ad_{D_u}$ describes a $\Zee/2$-grading of $\fder\, \cO(1;\underline{n})$.

The linear map $\ad_{D_u}$ sends elements of $\fvect(1;\underline{n})$ to $\fvect^{(1)}(1;\underline{n})$ because $D_u$ is a linear combination of a vector field and elements of the form $\del^{2^i}$, and
\begin{gather*}
\big[\del^{2^i}, f\del\big]= \big(\del^{2^i}f\big)\del = \big[\del, (\del^{2^i-1}f)\del\big]\in\fvect^{(1)}(1;\underline{n}) \qquad \text{for any $f\in \cO(1;\underline{n})$.}
\end{gather*}
This means that the restriction $(\ad_{D_u})|_{\fvect(1;\underline{n})}$ describes a $\Zee/2$-grading of $\fvect(1;\underline{n})$ and the restriction $(\ad_{D_u})|_{\fvect^{(1)}(1;\underline{n})}$ describes a~$\Zee/2$-grading of $\fvect^{(1)}(1;\underline{n})$. In what follows, we will denote the latter restriction by $D_u$ as well. Different polynomials $u$ correspond to different gradings, as $[D_u, x\del] = u\del$. We need, however, not individual gradings, but their equivalence classes. We were unable to solve this problem completely so far.

In notation of~\cite{T}, any torus of the restricted closure of $\fvect(1;\underline{n})$ lying in the maximal subalgebra of elements of non-negative degree (assuming $\deg x=1$) is called an \textit{inner} one, the other tori are called \textit{outer} ones.

If $u(0)=a\neq0$, the derivation~$D_u$ is an outer toroidal derivation of $\cO(1;\underline{n})$ (i.e., it spans an outer torus of $2$-closure of $\fvect(1;\underline{n})$). The automorphism group of $\fvect^{(1)}(1;\underline{n})$ acts transitively on the set of outer toroidal derivations, see~\cite{T}.\footnote{\label{footDerDer}In~\cite{T}, Tyurin proved that, for any two outer toroidal derivations~$O_1$ and~$O_2$ of $\fvect(1;\underline{n})$, there exists an $A\in \operatorname{Aut}\big(\fvect(1;\underline{n})\big)$ such that $O_2 = A^{-1}O_1A$. Instead, we need ``for any two outer toroidal derivations~$O_1$ and~$O_2$ of $\fvect^{(1)}(1;\underline{n})$, there exists an $A\in \operatorname{Aut}\big(\fvect^{(1)}(1;\underline{n})\big)$ such that $O_2 = A^{-1}O_1A$''. This is so because ``our'' $D_u = \big(D^T_u\big)|_{\fvect^{(1)}}\in\fder\, \fvect^{(1)}(1;\underline{n})$ is the restriction of the corresponding ``Tyurin's'' $D^T_u\in\fder\, \fvect(1;\underline{n})$ to $\fder\,\fvect^{(1)}(1;\underline{n})$.

If $D_u$ and $D_v$ are outer toroidal derivations of $\fvect^{(1)}(1;\underline{n})$, then $D^T_u$ and $D^T_v$ are outer toroidal derivations of $\fvect(1;\underline{n})$, and therefore there exists an $A\in \operatorname{Aut}(\fvect(1;\underline{n}))$ such that $D^T_v = A^{-1} D^T_u A$. Any automorphism of $\fvect(1;\underline{n})$ preserves $\fvect^{(1)}(1;\underline{n})$, and hence ${D_v = (\bar A)^{-1} D_u \bar A}$, where $\bar A:=A|_{\fvect^{(1)}}(1;\underline{n})$ is an automorphism of $\fvect^{(1)}(1;\underline{n})$.} Hence, any derivation~$D_u$ with $a\neq 0$ is conjugate by an automorphism of $\fvect^{(1)}(1;\underline{n})$ to the derivation
\begin{gather*}
{D_1 = (1+x)\del + \sum\limits_{1\leq i\leq n-1}\del^{2^i}}.
\end{gather*}
Recall that
$
(u+v)^{(k)} \!=\! \sum\limits_{0\leq i\leq k}\!u^{(k-i)}v^{(i)}
$
and for $k\!=\!-2, \dots, 2^{n-1}-2$, define $e_k\!\in\! \big(\fvect^{(1)}(1;\underline{n}); D_1\big)^{\langle 1\rangle}$ by setting
\begin{gather}
e_{-2} = \del^2 +(1+x)\del,\nonumber\\
e_{k}= \del\big( \big(x+x^{(2)}\big)^{(k+2)}\big)\del=(1+x)\big(x+x^{(2)}\big)^{(k+1)}\del \qquad \text{for $k>-2$},\label{evPart}
\end{gather}
and for $k=-1,\ldots,2^{n-1}-2$, define $o_k$ by formulas
\begin{gather}\label{odPart}
o_{k} = \big(x+x^{(2)}\big)^{(k+1)}\del, \qquad \text{where} \ \ k=-1,\dots,2^{n-1}-2.
\end{gather}
Direct computations show that $D_1(e_k) = 0$ and $D_1(o_m) = o_m$. Let $e_{-2},\ldots,e_{2^{n-1}-2}$ and $\Pi o_{-1},\ldots,\Pi o_{2^n-2}$ form a basis in the Lie superalgebra corresponding to the $\Zee/2$-grading
defined by $D_1$. The formulas{\samepage
\begin{gather*}\begin{split}&
 e_k \mapsto X_k,\qquad \text{where} \ \ k=-2,\ldots,2^{n-1}-2,\\
 & o_m \mapsto Y_m,\qquad \text{where} \ \ m=-1,\ldots,2^{n-1}-2,\end{split}
\end{gather*}
establish an isomorphism between this Lie superalgebra, spanned by the~$e_i$ and~$o_j$, and~$\fkl_{n-1}$.}

If $u(0)=a=0$, the derivation~$D_u$ spans an inner torus. Hence, see~\cite{T}, the derivation~$D_u$ is conjugate to the derivation $D_f$, where $f = \sum\limits_{1\leq i\leq n-1} c_ix^{(2^i)}$ and the summands $c_ix^{(2^i)}$ are the corresponding terms of $u$. This gives us an $(n-1)$-parametric family of $\Zee/2$-gradings. Let the automorphism~$\sigma_{\veps}$ of $\cO(1;\underline{n})$ be given on its generators by the formulas
\begin{gather*}
 \sigma_\veps(x) = \veps x,\quad \sigma_\veps\big(x^{(2)}\big) = \veps^2 x^{(2)},\quad \ldots,\quad \sigma_\veps \big(x^{(2^{n-1})}\big) = \veps^{2^{n-1}} x^{(2^{n-1})},
\end{gather*}
where $\veps\in{\mathbb K}^{\times}$. Then
\begin{gather*}
\sigma_\veps\bigg(\bigg(\sum\limits_{1\leq i\leq n-1} c_{i}x^{(2^i)}\bigg)\partial\bigg)=\bigg(\sum\limits_{1\leq i\leq n-1} c_{i}\veps^{2^i-1}x^{(2^i)}\bigg)\partial.
\end{gather*}
Hence, the conjugacy class of the derivation~$D_u$ is defined by a tuple of parameters~$(c_1,c_2,\ldots,$ $c_{n-1})$ up to an equivalence of the form
\begin{equation}\label{rescaleParm}
(c_1,c_2,\ldots,c_{n-1})\sim\big(\veps c_1, \veps^3 c_2, \ldots, \veps^{2^{n-1}-1}c_{n-1}\big).
\end{equation}
Therefore, the number of parameters of $\Zee/2$-gradings reduces to $n-2$. In particular, for $n=2$, we do not, indeed, have parametric families of $\Zee/2$-gradings, see Section~\ref{ssW2}.

Observe that $x\cdot \del u = 0$, so
\begin{gather*}
 0 = \del( x \del u ) = \del u + x \del^2 u, \qquad \text{which implies} \quad x\del^2 u = \del u.
\end{gather*}
Then
\begin{gather*}
D_u=\big(u+x+xu\del^2 u\big)\del = (u+x+u\del u)\del.
\end{gather*}
Desuperization of any superization of $\fvect^{(1)}(1;\underline{n})$ can be considered as a Lie subalgebra in $\fv_{n+1}:=\big(\fvect^{(1)}(1;\underline{n})\big)^{\langle 1\rangle}$, see~\eqref{fv_n}. The operator $\ad_{D_u}$ preserves $\fv_{n+1}$, so the grading it defines on $\fvect^{(1)}(1;\underline{n})$ can be extended to $\fv_{n+1}$ with the same definition: $D\mapsto [D_u,D]$. This extended operator describes a $\Zee/2$-grading of $\fv_{n+1}$. The elements of $\fv_{n+1}$ are said to be \textit{$u$-even} or \textit{$u$-odd} if they are even or odd, respectively, in this grading.

In particular, for the grading on $\fv_{n+1}$ given by the function $u=0$, we have $D_0=x\del$, and the resulting superization is isomorphic to $\fk(1; \underline{n-1}|1)$, see~\cite[Section~7]{BLLS1}.

Consider the linear maps $T_u, A_u\colon \fv_{n+1} \tto \fv_{n+1} $ defined as follows (their action on other elements being extended by linearity):
\begin{gather*}
T_u (f\del) = (f + u(1+\del u)\del f)\big(1+\del u + u\del^2 u\big)\del \qquad \text{for any} \ f\in\cO(1;\underline{n}), \\
T_u \del^2 = (T_u \del)^2 = \del^2 + \big(\del^2 u + u\del^3 u + u\big(\del^2 u\big)^2 + u\del u \big(\del^2 u\big)^2\big)\del, \\
A_u (f\del) = (f + u\del u\del f)\del \qquad \text{for any} \ f\in\cO(1;\underline{n}), \qquad
A_u \del^2 = \del^2.
\end{gather*}
For any $X\in \fv_{n+1} $, one can check that $T_u X$ is $u$-even (resp.~$u$-odd) if and only if $X$ is $0$-even (resp.~$0$-odd). We omit the calculations that show it, but the idea is as follows: $\del$ is $0$-odd, while $\big(1+\del u + u\del^2 u\big)\del$ is $u$-odd. Besides, if we extend the concept of $u$-evenness/$u$-oddness to $\cO(1;\underline{n})$, then $f + u(1+\del u)\del f$ is $u$-even (resp. $u$-odd) for any $f\in\cO(1;\underline{n})$ if and only if $f$ is $0$-even (resp.~$0$-odd).

Also, the following is true for any $X,Y\in \fvect(1;\underline{n})$:
\begin{gather*}\begin{split}
 & [T_u X, T_u Y] = T_u A_u [X,Y],\qquad
 (T_u X)^2 = T_u A_u X^2 \qquad \text{for any $0$-odd} \ X,\\
& \big[T_u\del^2, T_u X\big] = T_u\big(\big(1+u\del^2 u\big)\cdot \big[\del^2, X\big]\big).\end{split}
\end{gather*}

Note also that the operator $T_u$ is invertible, so its image is the whole $\fv_{n+1} $. The invertibility of~$T_u$ follows from the fact that the matrix of~$T_u$ is upper triangular with $1$'s on the main diagonal in the basis $\del^2, \del, x\del, \dots, x^{(2^n-1)}\del$.

This means that the superization given by any function $u$ such that $u(0)=0$ can be considered as a deform of the superization given by $u=0$ with the deformation parametrized by the polynomial~$u$ or, which is the same, the coefficients of $u$ as follows:
\begin{gather}
{}[X,Y]_u = \begin{cases} A_u[X,Y] & \text{if} \ X,Y\in \fvect(1;\underline{n}),\\
\big(1+u\del^2 u\big)\cdot [X,Y] & \text{if} \ X = \del^2 \ \text{and} \ Y\in \fvect(1;\underline{n}), \\
\text{defined by linearity and anti-symmetry} & \text{in other cases},\end{cases}\nonumber\\
\big(X^{2}\big)_u = A_u X^{2} \qquad \text{for 0-odd $X$}.\label{defBracket}
\end{gather}

This deformation is filtered relative the decreasing filtration in which $\cL_{-2}=\fs(\fv_{n+1})$, see~\eqref{S(g)}, while $\cL_{k-1}$ for $k\geq -1$ consists of vector fields of the form $f\del$, where $f\in \cO(1;\underline{n})$ does not contain term of degree $<k$. Such superizations are listed in heading 2) of Theorem~\ref{T_vect}.

Finally, observe that these formulas do not capture the trivial $\Zee/2$-grading given by the derivation $U = 0$; the even part~$\fg_{\ev}$ of this grading is the whole Lie algebra~$\fvect^{(1)}(1;\underline{n})$. Since any torus in $\fder \big(\fvect^{(1)}(1;\underline{n})\big)$ is of the form~$D_u$, see~\cite{Kuz,T} and footnote~\ref{footDerDer} on p.~\pageref{footDerDer}, we completely described all $\Zee/2$-gradings of~$\fvect^{(1)}(1;\underline{n})$ and the corresponding superizations.
\end{proof}

\begin{Remark} Observe an unpredicted fact: according to equations~\eqref{defBracket}, if $X,Y\in \fvect(1;\underline{n})$, then $[X,Y]_u$ can be expressed in terms of $[X,Y]$ and $u$.
\end{Remark}

In particular, concerning the $0$-even part of $\fvect(1;\underline{n})$, this implies the following fact:
\begin{Corollary}\label{Cor_Solvable}If $u(0)=0$, then the even part of the superization corresponding to~$D_u$ is a~solvable Lie algebra.\footnote{For more examples of simple Lie superalgebras whose even parts are solvable, see~\cite{BGL}. This is a phenomenon indigenous to $p=2$.}
\end{Corollary}

For $n=2$, any outer torus is conjugate to an inner one, see \eqref{*}. This is not so for $n>2$.

\begin{Lemma}\label{sssCj1} Let $n>2$, consider the superization $\fkl_{n-1}=\fg_{\ev}\oplus\fg_{\od}$ corresponding to the outer
torus~$D_1$. Its even part~$\fg_{\ev}\simeq\fv_{n}$ is spanned by the elements \eqref{evPart}, and
~$\fg_{\od}$ is spanned by the elements \eqref{odPart}. The odd part~$\fg_{\od}$ is a reducible
$\fg_{\ev}$-module with no lowest weight vector and with the two highest weight vectors~$\Pi o_{2^{n-1}-3}$
and~$\Pi o_{2^{n-1}-2}$ with respect to the standard $\Zee$-grading of~$\fg_{\ev}$, namely
$(\fg_{\ev})_{k} = {\mathbb K} e_k$ for $k=-2, \dots, 2^{n-1}-2$.
\end{Lemma}
\begin{proof}Using~\eqref{evPart} and~\eqref{odPart}, for positive generators $e_{2^k-1}$, where $k=1,\ldots,n-2$, of~$\fg_{\ev}$ we have
\begin{gather*}
 \big[e_{2^k-1}, \Pi o_{2^{n-1}-2}\big] = \Pi \big[ (\del w) w^{(2^k)}\del, w^{(2^{n-1}-1)}\del\big] \\
\hphantom{\big[e_{2^k-1}, \Pi o_{2^{n-1}-2}\big]}{} = \Pi \big( w^{(2^k)}w^{(2^{n-1}-2)} + w^{(2^k)}w^{(2^{n-1}-1)} + w^{(2^k-1)} w^{(2^{n-1}-1)} \big) \del,
\end{gather*}
where $w = \big(x+x^{(2)}\big)$. Note that
 \begin{gather*}
 2^{n-1}-1 = 1\cdot 2^{n-2} + 1\cdot 2^{n-3} + \dots + 1\cdot 2 + 1,\\
 2^{n-1}-2 = 1\cdot 2^{n-2} + 1\cdot 2^{n-3} + \dots + 1\cdot 2 + 0.
 \end{gather*}
 Then by Lucas theorem for $k<n-1$, we have the following equalities, where the marked factor occupies $(k+1)$st position from the right in the first two lines, while it stands in the $k$th position from the right in the third line, where the mark is under the place at which the 1s in the bottom of the binomial coefficients start to appear
 \begin{gather*}
 \binom{2^{n-1} -2 + 2^k}{2^k} \mod 2 = \binom{1}{0} \binom{0}{0} \cdots \binom{0}{0} \underbrace{\binom{0}{1}} \binom{1}{0} \cdots \binom{1}{0}\binom{0}{0} = 0,\\
 \binom{2^{n-1}-1 + 2^k}{2^k} \mod 2 =\binom{1}{0} \binom{0}{0} \cdots \binom{0}{0} \underbrace{\binom{0}{1}} \binom{1}{0} \cdots \binom{1}{0}\binom{1}{0} = 0,\\
 \binom{2^{n-1}-1 + 2^k - 1}{2^k - 1} \mod 2 =\binom{1}{0} \binom{0}{0} \cdots \binom{0}{0}\underbrace{\binom{1}{1}} \cdots \binom{1}{1}\binom{0}{1} = 0.
 \end{gather*}
Thus, we have $w^{(2^k)}w^{(2^{n-1}-2)} =0$, $w^{(2^k)}w^{(2^{n-1}-1)}=0$, and $w^{(2^k-1)} w^{(2^{n-1}-1)} = 0$. Finally, we see that $[e_{2^k-1}, \Pi o_{2^{n-1}-2}] = 0$ for $k=1,\ldots,n-2$. Computations of the same kind show that $[e_{2^k-1}, \Pi o_{2^{n-1}-3}] = 0$ for any $k=1,\ldots,n-2$. Therefore, $o_{2^{n-1}-2}$ and $o_{2^{n-1}-3}$ are highest weight vectors in~$\fg_{\od}$.

Now, let us prove that there are no lowest weight vectors in~$\fg_{\od}$. Consider the action of $e_{-1}$ on elements of $\fg_{\od}$. Using~\eqref{evPart} and~\eqref{odPart}, we obtain
\begin{gather}\label{outerDelModule}
[e_{-1},\Pi o_{-1}] = o_1,\qquad [e_{-1}, \Pi o_k] = \Pi o_{k} + \Pi o_{k-1} \qquad \text{for any $k=0,\ldots,2^{n-1}-2$}.
\end{gather}
Therefore, the equation
\begin{gather*}{}
\bigg[e_{-1}, \sum \limits_{-1\leq k\leq 2^{n-1}-2} q_k\Pi o_k\bigg] = 0
\end{gather*}
for lowest weight vectors reduces to the following system of linear equations
 \begin{gather*}
 q_{-1} + q_0 = 0, \quad q_0 + q_1 = 0, \quad \ldots, \quad q_{2^{n-1}-3} + q_{2^{n-1}-2} = 0, \quad q_{2^{n-1}-2} =0.
 \end{gather*}
Clearly, this system has only one solution: $q_i= 0$ for $i=-1,\ldots,2^{n-1}-2$. Hence, there are no lowest weight vectors in~$\fg_{\od}$.

Observe that the highest weight vector $\Pi o_{2^{n-1}-2}$ generates the whole $\fg_{\ev}$-module~$\fg_{\od}$ and the highest weight vector $\Pi o_{2^{n-1}-3}$ generates a $\fg_{\ev}$-submodule of codimesion~1 in~$\fg_{\od}$, as follows from~\eqref{outerDelModule}.
\end{proof}

In what follows we consider only generating functions of the form
\begin{gather}
u(c) = c_0 + \sum\limits_{1\leq k\leq n-1} c_kx^{(2^k)},\label{param}
\end{gather}
where $c = (c_0, c_1, \ldots, c_{n-1})$ is a set of free parameters, cf.~\eqref{rescaleParm}.
In what follows by ``grading $c=(c_0, c_1, \ldots, c_{n-1})$'' we mean the grading given by the corresponding $D_{u(c)}$. One should bear in mind that we call these parameters \textit{free} by an abuse of the language, since they are classes modulo equivalence \eqref{rescaleParm}.

Let us consider the classification of $\Zee/2$-gradings of~$\fvect^{(1)}(1;\underline{2})$ in more detail.

\subsection[Superizations of $\mathfrak{g}=\mathfrak{vect}^{(1)}(1;\underline{2})$]{Superizations of $\boldsymbol{\mathfrak{g}=\mathfrak{vect}^{(1)}(1;\underline{2})}$}\label{ssW2}
The general solution~$U_{\text{lin}}$ of linear equations~\eqref{cond2lin}, where $c_{1,1}$, $c_{1,2}$, $c_{1,3}$, $c_{2,1}$, $c_{3,1}$ are free parameters, and the general solution~$U$ of equation~\eqref{cond2}, where $c_0$,
$c_1$ are free parameters, and the schematic form of $U$ are as follows:{\samepage
\begin{gather*}
U_{{\rm lin}} =
\begin{pmatrix}
 c_{1,1} & c_{1,2} & c_{1,3} \\[-1pt]
 c_{2,1} & 0 & c_{1,2} \\[-1pt]
 c_{3,1} & c_{2,1} & c_{1,1}
\end{pmatrix},
\qquad U = \left(
\begin{matrix}
 c_0 c_1+1 & c_0 & c_0^2 \\[-1pt]
 c_1 & 0 & c_0 \\[-1pt]
 c_1^2 & c_1 & c_0 c_1+1
\end{matrix}
\right),\qquad \begin{pmatrix}
 * & * & * \\[-1pt]
 * & & * \\[-1pt]
 * & * & *
\end{pmatrix}.
\end{gather*}
This solution corresponds to the derivation~$D_u$ with $u(c)=c_0 + c_1x^{(2)}$.}

{\bf Parameters and $\boldsymbol{\Aut\fg}$.} We denote by (see \eqref{param})
\begin{gather}\label{paramRec}
\fg_{\ev}(c):=\Ker D_{u(c)}
\end{gather}
the even part of the corresponding $\Zee/2$-grading. For example, for $\fg=\fvect^{(1)}(1;\underline{3})$, by $\fg_{\ev}(c_0,c_1)$ we mean the even part of the $\Zee/2$-grading $D_u$ where $u=c_0+ c_1x^{(2)}$.

Certain automorphisms and their actions:
\begin{gather}
\begin{pmatrix}
 \frac{1}{c_0} & 0 & 0 \\
 0 & 1 & 0 \\
 \frac{c_1 + c_0 c_1}{c_0} & 0 & c_0
 \end{pmatrix}
\qquad \text{maps $D_{u(c_0,c_1)}$ to $D_{u(1,c_1)}$ for all $c_0\neq 0$},\nonumber\\
\begin{pmatrix}
 1+c_1 & 0 & 1\\
 0 & 1 & 0 \\
 1 & 0 & 0
 \end{pmatrix}\qquad \text{maps $D_{u(1,c_1)}$ to $D_{u(1,1)}$ for all $c_1\in{\mathbb K}$},\nonumber\\
\begin{pmatrix}
 c_1 & 0 & 0 \\ 0 & 1 & 0 \\ 0 & 0 & \frac{1}{c_1}
 \end{pmatrix}\qquad \text{maps $D_{u(0,c_1)}$ to $D_{u(0,1)}$ for all $c_1\neq0$},\nonumber\\
\label{*}
\begin{pmatrix}
 1 & 0 & 1 \\ 0 & 1 & 0 \\ 0 & 0 & 1
 \end{pmatrix}\qquad \text{maps $D_{u(1, 1)}$ to~$D_{u(0,1)}$}.
 \end{gather}
Observe that the automorphism \eqref{*} sends the outer derivation $D_{1+x^{(2)}}$ to the inner derivation~$D_{x^{(2)}}$.

The final answer: there are \textit{three} inequivalent $\Zee/2$-gradings of $\fvect^{(1)}(1;\underline{2})\simeq{\mathfrak{o}}^{(1)}(3)$; their $\fg_{\ev}$'s are
\begin{gather}
\fg_{\ev}(0,0)={\mathbb K} \eb_0 \qquad \text{and}\qquad \fg_{\ev}(1,1) ={\mathbb K} (\eb_{-1} + \eb_0 + \eb_1),\label{three}
\end{gather}
and the trivial $\Zee/2$-grading with $\fg_{\ev}=\fvect^{(1)}(1;\underline{2})$, not \textit{two} as is the case for $\Zee/2$-gradings of ${\mathfrak{o}}^{(1)}(3)={\mathfrak{o}}(3)$ for any $p\neq2$.

{\bf Lie superalgebras corresponding to the gradings (\ref{three}).}

{\it Grading $c=(0,0)$.} The minimal Lie subsuperalgebra
\begin{gather}\label{subSup}
\fg^{\langle 1\rangle} := \big(\fg^{\langle 1\rangle}\big)_{\ev}\oplus\big(\fg^{\langle 1\rangle}\big)_{\od}
\end{gather}
of the restricted closure $\bar\fg$ containing all elements $v^{[2]}$, where
$v\in\fg_{\od}=\big(\fg^{\langle 1\rangle}\big)_{\od}$, see equation~\eqref{resCl}, is as follows:
\begin{gather*}
 \big(\fg^{\langle 1\rangle}\big)_{\ev} = \Span\big(e_{-2} = \del^2,\, e_0 = x \del,\, e_2 = x^{(3)}\del\big),\\
 \big(\fg^{\langle 1\rangle}\big)_{\od} = \Span\big(o_{-1} = \del,\, o_1 = x^{(2)} \del\big).
\end{gather*}
The squaring in $\fs\big(\fg^{\langle 1\rangle}\big)$ is given by the following formulas:
\begin{gather*}
 o_{-1}^2 = (\del)^2 = e_{-2} = \del^2, \qquad o_1^2 = \big(x^{(2)}\del\big)^2 = e_2 = x^{(3)}\del.
\end{gather*}
The commutation relations in~$\big(\fg^{\langle 1\rangle}\big)_{\ev}$ are defined by the following formulas:
\begin{gather*}
 [e_0,e_{-2}] = 0,\qquad [e_0,e_2]= 0,\qquad [e_{-2},e_2] = e_0.
\end{gather*}
Observe that~$\big(\fg^{\langle 1\rangle}\big)_{\ev}$ is the Heisenberg algebra. The $\big(\fg^{\langle 1\rangle}\big)_{\ev}$-module structure in~$\big(\fg^{\langle 1\rangle}\big)_{\od}$ is given by the following formulas:
\begin{alignat*}{4}
 &[e_{-2}, o_{-1}] = 0, \qquad&& [e_0, o_{-1}] = o_{-1}, \qquad&& [e_2, o_{-1}] = o_1,&\\
 & [e_{-2}, o_1] = o_{-1}, \qquad&& [e_0, o_1] = o_1, \qquad&& [e_2, o_1] = 0.&
\end{alignat*}

It is easy to see that~$\fs\big(\fg^{\langle 1\rangle}\big)\simeq{\mathfrak{oo}}^{(1)}_{I\Pi}(1|2)\simeq\fk(1;\underline{1}|1)$; for the (non-obvious) definition of the~$\fk(1;\underline{n}|1)$, see~\cite{BGLLS2}. The Lie
superalgebra~${\mathfrak{oo}}^{(1)}_{I\Pi}(1|2)$ is given by supermatrices of the form
\begin{gather*}
 \begin{pmatrix}
 0 & b_2 & b_1 \\
 b_1 & a_1 & a_2 \\
 b_2 & a_3 & a_1
 \end{pmatrix}, \qquad \text{where} \ \ a_i,b_i\in{\mathbb K}.
\end{gather*}
The correspondence between the abstract and matrix representations of the elements of~${\mathfrak{oo}}^{(1)}_{I\Pi}(1|2)$ is as follows
\begin{gather*}
 e_{-2} = e_{2,3}, \qquad e_0 = e_{2,2}+e_{3,3},\qquad e_2 = e_{3,2},\qquad
 o_{-1} = e_{1,3} + e_{2,1},\qquad o_1 = e_{1,2} + e_{3,1}.
\end{gather*}

{\it Grading $c=(1,1)\simeq(1,0)$.} The Lie subsuperalgebra $\fg^{\langle 1\rangle}$, see equation~\eqref{subSup}, is as follows
\begin{gather*}
 \big(\fg^{\langle 1\rangle}\big)_{\ev} = \Span\big(e_1 = \big(1+x+x^{(2)}\big)\del,\, e_2 = \del^2 + (1+x)\del,\, e_3 = \del^2 + \big(x+x^{(3)}\big)\del\big),\\
 \big(\fg^{\langle 1\rangle}\big)_{\od} = \Span\big(o_1 = (1+x)\del,\, o_2 = \big(1+x^{(2)}\big)\del\big).
\end{gather*}
The squaring in $\fs\big(\fg^{\langle 1\rangle}\big)$ is given by the following formulas:
\begin{gather*}
 o_1^2 = e_2, \qquad o_2^2 = e_3.
\end{gather*}
The commutation relations in~$\big(\fg^{\langle 1\rangle}\big)_{\ev}$ are defined by the following formulas:
\begin{gather}\label{neHei}
 [e_1, e_2] = e_1,\qquad [e_1,e_3] = 0,\qquad [e_2,e_3] = e_1.
\end{gather}
The $\big(\fg^{\langle 1\rangle}\big)_{\ev}$-module structure in~$\big(\fg^{\langle 1\rangle}\big)_{\od}$ is given by the following formulas:
\begin{alignat*}{4}
& [e_1, o_1] = o_1+o_2, \qquad && [e_2, o_1] = 0, \qquad&& [e_3, o_1] = o_2,&\\
& [e_1, o_2] = o_2, \qquad&& [e_2, o_2] = o_1+o_2, \qquad && [e_3, o_2] = 0.&
\end{alignat*}
Let us show that $\fs\big(\fg^{\langle 1\rangle}\big)\simeq{\mathfrak{oo}}^{(1)}_{II}(1|2)$. The Lie superalgebra~${\mathfrak{oo}}^{(1)}_{II}(1|2)$ consists of symmetric supermatrices with supertrace~0. For a basis of
its even part we take
\begin{gather*}
E_1 = e_{1,1}+e_{2,2},\qquad E_2 = e_{2,2}+e_{3,3},\qquad E_3=e_{2,3}+e_{3,2}.
\end{gather*}
For a basis of its odd part we take $O_1 = e_{1,2} + e_{2,1}$, $O_2=e_{1,3}+e_{3,1}$. The commutation relations are given by the following formulas:
\begin{gather*}
 [E_1, E_2] = 0,\qquad [E_1, E_3] = E_3,\qquad [E_2, E_3] = 0.
\end{gather*}
The squaring in $\fs\big(\fg^{\langle 1\rangle}\big)$ is given by the following formulas:
\begin{gather*}
 O_1^2 = E_1,\qquad O_2^2 = E_1 + E_2.
\end{gather*}
The $\big(\fg^{\langle 1\rangle}\big)_{\ev}$-module structure is as follows
\begin{gather*}
 [E_2, O_1] = O_1, \qquad [E_2, O_2] = O_2.
\end{gather*}
The multiplication by elements~$E_1$ and~$E_3$ is given by the following formulas:
\begin{gather*}
 [E_1, O_1] = 0,\qquad [E_1, O_2] = O_2,\qquad [E_3, O_1] = O_2,\qquad [E_3, O_2] = O_1.
\end{gather*}
The isomorphism between~${\mathfrak{oo}}^{(1)}_{II}(1|2)$ and~$\fg^{\langle 1\rangle}$ is given by the following formulas:
\begin{gather*}
 e_1 = E_3,\qquad e_2 = E_1,\qquad e_3 = E_2 + E_3,\qquad o_1 = O_1,\qquad o_2 = O_1 + O_2.
\end{gather*}

\begin{Remark}\label{R_oII} \quad
\begin{enumerate}\itemsep=0pt
\item[1)] Observe that the even parts of both superizations are solvable Lie algebras, but the corresponding Lie superalgebras are simple. For more details and examples of this phenomenon indigenous to $p=2$, see Shchepochkina's comment (the last section) in~\cite{BGL}.
\item[2)] Observe that due to Theorem~\ref{T_vect} the Lie superalgebra~${\mathfrak{oo}}_{II}^{(1)}(1|2)$ is a filtered deform of the Lie superalgebra~${\mathfrak{oo}}^{(1)}_{I\Pi}(1|2)$.
\end{enumerate}
\end{Remark}

\subsection{Superizations corresponding to inner tori} Consider a superization of $\fvect^{(1)}(1;\underline{n})$ given by a generating function $u$ such that $u(0)=0$. For such a function $u$, the derivation~$D_u$ spans an inner torus. According to~\cite{T}, the derivation~$D_u$ is conjugate to the derivation $D_f$ in $\fvect^{(1)}(1;\underline{n})$, where $f = \sum\limits_{1\leq i\leq n-1} c_ix^{(2^i)}$ and the $c_ix^{(2^i)}$ are the corresponding terms of $u$, so the superizations corresponding to $D_u$ and $D_f$ are isomorphic. For this reason, in the rest of this section we consider only generating functions of the form $\sum\limits_{1\leq i\leq n-1} c_ix^{(2^i)}$. This gives us an $(n-1)$-parametric family of $\Zee/2$-gradings.

\begin{Conjecture}\label{inner-tori-conj} Let $u_1=\sum\limits_{1\leq k\leq n-1} c_k x^{(2^k)}$ and $u_2= \sum\limits_{1\leq k\leq n-1} b_k x^{(2^k)}$ be generating functions, then the superizations corresponding to $D_{u_1}$ and $D_{u_2}$ are isomorphic if and only if there exists $\veps\in{\mathbb K}^\times$ such that $c_k = \veps^{2^k-1} b_k$ for all $k=1,\dots, n-1$.
\end{Conjecture}

The following is a sketch of a proof of this conjecture, including a complete computer-aided proof for $2\leq n\leq 6$.

First of all, if there is an $\veps\neq0$ such that $c_k = \veps^{2^k-1} b_k$ for all $k=1,\dots, n-1$, then $D_{u_1}$ and $D_{u_2}$ are conjugate by the automorphism of $\fvect^{(1)}(1;\underline{n})$ given by $x^{(m)}\del \mapsto \veps^{m-1}x^{(m)}\del$ (this automorphism is generated by the automorphism of $\cO(1;\underline{n})$ given by $x^{(m)} \mapsto \veps^{m}x^{(m)}$), and therefore the corresponding superizations are isomorphic.

By Theorem~\ref{T_vect}, the superization given by $u_1$ is isomorphic to a deform of $\fk(1; \underline{n-1}|1)$. Consider the even part of such a deform. According to~\eqref{defBracket}, it contains a commutative subalgebra of codimension $1$, which is the $0$-even part of $\fvect(1;\underline{n})$. The element $\del^2$ acts on this commutative subalgebra, and its action is given by
 \begin{gather}\label{inner-torui-d2-action}
 v\del \mapsto \big(1+u\del^2 u\big)\del^2v\del \qquad \text{for any $0$-odd $v\in \cO(1;\underline{n})$}.
 \end{gather}

\begin{Conjecture}[proved by computer for $n=2, 3, 4, 5, 6$]\label{char-pol} The characteristic polynomial of the linear operator \eqref{inner-torui-d2-action} on the $0$-even part of $\fvect(1;\underline{n})$ for $n\geq 2$ is defined by the formula
 \begin{gather*}
 \lambda^{2^{n-1}} + \sum\limits_{0\leq k\leq n-2} c_{n-1-k}^{2^{k+1}} \lambda^{2^k}.
 \end{gather*}
\end{Conjecture}

For two deforms given by $u_1$ and $u_2$ to be isomorphic to each other, their even parts have to be isomorphic, which means that the two actions of $\del^2$ on the $0$-even part of $\fvect(1;\underline{n})$ have to be conjugate (similar) up to a non-zero scalar factor, i.e., if $A$ and $A'$ are such operators, then $A' = \alpha MAM^{-1}$ for some non-zero $\alpha\in{\mathbb K}$ and an invertible linear map $M$. If two operators are conjugate up to a scalar multiple $\alpha$, the roots of their characteristic polynomials differ by the same multiple, so if the dimension of the space they act on is $d$ (in our case, $d = 2^{n-1}$) and one polynomial is equal to $\sum\limits_{0\leq i\leq d} z_i\lambda^i$, where $z_i\in{\mathbb K}$ and $z_d=1$, then the other polynomial would have the form $\sum\limits_{0\leq i\leq d} z_i\alpha^{d-i}\lambda^i$. So, if Conjecture \ref{char-pol} is correct, then for the superizations given by $u_1$ and $u_2$ to be isomorphic, there must exist a non-zero $\alpha\in{\mathbb K}$ such that
\begin{gather*}
\lambda^{2^{n-1}} + \sum\limits_{0\leq k\leq n-2} c_{n-1-k}^{2^{k+1}} \lambda^{2^k} = \lambda^{2^{n-1}} + \sum\limits_{0\leq k\leq n-2} \alpha^{2^{n-1}-2^k}b_{n-1-k}^{2^{k+1}} \lambda^{2^k},
\end{gather*}
or, equivalently, $c_k = \veps^{2^k-1} b_k$, where $\veps = \sqrt{\alpha}$, for all $k=1,\dots, n-1$. Thus, Conjecture~\ref{inner-tori-conj} follows from Conjecture~\ref{char-pol}.

So, assuming that Conjecture \ref{char-pol} is true, the set of equivalence classes of superizations of $\fvect^{(1)}(1;\underline{n})$ corresponding to generating functions $u$ such that $u(0)=0$ consist of the following two types:
\begin{enumerate}\itemsep=0pt
\item[A)] the superization corresponding to $u=0$, which is isomorphic to $\fk(1; \underline{n-1}|1)$,
\item[B)] an $(n-2)$-parametric family of its pairwise non-isomorphic deforms. Note, though, that it is not a result of $(n-2)$-parametric deformation of $\fk(1; \underline{n-1}|1)$; to obtain all these deforms, an $(n-1)$-parametric deformation is needed, but the deforms obtained from some sets of parameters are isomorphic: parameters $(c_1, \dots, c_{n-1})$ and $(b_1, \dots, b_{n-1})$ produce isomorphic deformations if and only if there exists an $\veps\in{\mathbb K}^\times$ such that $c_k = \veps^{2^k-1}b_k$ for all $k\in 1,\dots,n-1 $.
\end{enumerate}

\subsection[Superizations of~$\fvect^{(1)}(1;\underline{3})$]{Superizations of~$\boldsymbol{\fvect^{(1)}(1;\underline{3})}$} Let $\fg=\fvect^{(1)}(1;\underline{3})$.

{\it Grading $c=(0,0,0)$.} Recall our notation \eqref{paramRec}. The Lie subsuperalgebra $\fg^{\langle 1\rangle}$, see equation~\eqref{subSup}, is as follows:
\begin{gather*}
 \big(\fg^{\langle 1\rangle}\big)_{\ev} = \Span\big(e_{-2} = \del^2,\, e_0 = x\del,\, e_2 = x^{(3)}\del,\, e_4 = x^{(5)}\del,\, e_6 = x^{(7)}\del\big),\\
 \big(\fg^{\langle 1\rangle}\big)_{\od} = \Span\big(o_{-1} = \del,\, o_1 = x^{(2)}\del,\, o_3 = x^{(4)}\del,\, o_5 = x^{(6)}\del\big).
\end{gather*}
The squaring in $\fs\big(\fg^{\langle 1\rangle}\big)$ is given by the following formulas:
\begin{gather*}
 o_{-1}^2 = e_{-2}, \qquad o_1^2 = e_2, \qquad o_3^2 = e_6, \qquad o_5^2 = 0.
\end{gather*}
The nonzero commutation relations in~$\big(\fg^{\langle 1\rangle}\big)_{\ev}$ are given by the following formulas:
\begin{gather*}
 [e_2,e_{-2}]= e_0,\qquad [e_4, e_{-2}]= e_2,\qquad [e_{-2}, e_6]= e_4.
\end{gather*}
The $\big(\fg^{\langle 1\rangle}\big)_{\ev}$-module structure in~$\big(\fg^{\langle 1\rangle}\big)_{\od}$ is given by the following formulas:
\begin{alignat*}{6}
 &[e_{-2},o_{-1}]= 0,\qquad &&[e_0,o_{-1}]= o_{-1},\qquad &&[e_2,o_{-1}]= o_1,\qquad &&[e_4,o_{-1}]= o_3,\qquad &&[e_6,o_{-1}]= o_5,& \\
& [e_{-2},o_1]= o_{-1},\qquad &&[e_0,o_1]= o_1,\qquad &&[e_2,o_1]= 0,\qquad &&[e_4,o_1]= o_5,\qquad &&[e_6,o_1]= 0,& \\
& [e_{-2},o_3]= o_1,\qquad &&[e_0,o_3]= o_3,\qquad &&[e_2,o_3]= o_5,\qquad &&[e_4,o_3]= 0,\qquad &&[e_6,o_3]= 0, &\\
& [e_{-2},o_5]= o_3,\qquad &&[e_0,o_5]= o_5,\qquad &&[e_2,o_5]= 0,\qquad && [e_4,o_5]= 0,\qquad && [e_6,o_5]= 0.&
\end{alignat*}
Let $L_k$ and $L^{(k)}$ be given by the formulas
\begin{gather}\label{L_k}
 L_0 = L^0=\big(\fg^{\langle 1\rangle}\big)_{\ev},\qquad L_{k} = [L_0, L_{k-1}], \qquad L^{(k)} = \big[L^{(k-1)}, L^{(k-1)}\big] \qquad \text{for $k>0$}.
\end{gather}
We have
\begin{gather*}
 L_1 = \Span(e_0,\, e_2,\, e_4),\qquad L_2 = \Span(e_0,\, e_2),\qquad L_3 = {\mathbb K} e_0,\qquad L_4 = 0.
\end{gather*}

{\it Grading $c=(0,1,0)$.} The Lie subsuperalgebra $\fg^{\langle 1\rangle}$, see equation~\eqref{subSup}, is as follows:
\begin{gather*}
 \big(\fg^{\langle 1\rangle}\big)_{\ev} = \Span\big(e_1 = (x + x^{(2)})\del,\, e_2 = x^{(3)}\del,\, e_3 = \big(x^{(5)} + x^{(6)}\big)\del,\\
\hphantom{\big(\fg^{\langle 1\rangle}\big)_{\ev} = \Span\big(}{} e_4 = \del^2 + (1 + x)\del,\, e_5 = x^{(7)}\del\big),\\
 \big(\fg^{\langle 1\rangle}\big)_{\od} = \Span\big(o_1 = (1 + x)\del,\, o_2 = x^{(2)}\del,\, o_3 = \big(x^{(4)} + x^{(5)}\big)\del,\, o_4 = x^{(6)}\del\big).
\end{gather*}
The squaring in $\fs\big(\fg^{\langle 1\rangle}\big)$ is given by the following formulas:
\begin{gather*}
 o_1^2 = e_4,\qquad o_2^2 = e_2, \qquad o_3^2 = e_5, \qquad o_4^2 = 0.
\end{gather*}
The non-zero commutation relations in~$\big(\fg^{\langle 1\rangle}\big)_{\ev}$ are given by the following formulas:
\begin{gather*}
 [e_1, e_4]=e_1,\qquad [e_2, e_4]=e_1,\qquad [e_3, e_4]=e_2+e_3,\qquad [e_4, e_5]=e_3.
\end{gather*}
The $\big(\fg^{\langle 1\rangle}\big)_{\ev}$-module structure in~$\big(\fg^{\langle 1\rangle}\big)_{\od}$ is given by the following formulas:
\begin{alignat*}{6}
&[e_1, o_1] = o_1+o_2, \quad && [e_2, o_1] = o_2,\quad && [e_3, o_1] = o_3+o_4,\quad && [e_4, o_1] = 0,\quad && [e_5, o_1] = o_4,&\\
& [e_1, o_2] = o_2,\quad && [e_2, o_2] = 0,\quad && [e_3, o_2] = o_4,\quad && [e_4, o_2] = o_1+o_2,\quad && [e_5, o_2] = 0,& \\
& [e_1, o_3] = o_3+o_4,\quad && [e_2, o_3] = o_4, \quad && [e_3, o_3] = 0,\quad && [e_4, o_3] = o_2,\quad && [e_5, o_3] = 0,& \\
& [e_1, o_4] = o_4,\quad && [e_2, o_4] = 0,\quad && [e_3, o_4] = 0, \quad && [e_4, o_4] = o_3+o_4,\quad && [e_5, o_4] = 0.&
\end{alignat*}

{\it Grading $c=(0,1,\beta)$, where $\beta\neq0$.} The Lie subsuperalgebra $\fg^{\langle 1\rangle}$, see equation~\eqref{subSup}, is as follows:
\begin{gather*}
 \big(\fg^{\langle 1\rangle}\big)_{\ev} = \Span\big(e_1 = \big(x + x^{(2)} + \beta x^{(4)}\big)\del,\, e_2 = (x^{(3)} + \beta x^{(5)})\del,\, e_3 = \big(x^{(5)} + x^{(6)}\big)\del,\\
\hphantom{\big(\fg^{\langle 1\rangle}\big)_{\ev} = \Span\big(}{} e_4 = \del^2 + \big(1 + x + \beta x^{(2)}\big)\del,\, e_5 = \big(x^{(3)} + \beta x^{(5)} + \beta^2 x^{(7)}\big)\del\big),\\
 \big(\fg^{\langle 1\rangle}\big)_{\od} = \Span\big(o_1 = \big(1 + x + \beta x^{(3)}\big)\del,\, o_2 = \big(x^{(2)} + \beta x^{(4)}\big)\del,\\
 \hphantom{\big(\fg^{\langle 1\rangle}\big)_{\od} = \Span\big(}{} o_3 = \big(x^{(4)} + x^{(5)}\big)\del,\, o_4 = x^{(6)}\del\big).
\end{gather*}
The squaring in $\fs\big(\fg^{\langle 1\rangle}\big)$ is given by the following formulas:
\begin{gather*}
 o_1^2 = e_4, \qquad o_2^2 = e_5, \qquad o_3^2 = \tfrac{1}{\beta^2}e_2 + \tfrac{1}{\beta^2}e_5, \qquad o_4^2 = 0.
\end{gather*}
The non-zero commutation relations in~$\big(\fg^{\langle 1\rangle}\big)_{\ev}$ are defined by the following formulas:
\begin{gather*}\begin{split}&
[e_1, e_4] = e_1 + \beta e_2, \\
& [e_2, e_4] = e_1 + \beta e_2 + \beta^2 e_3,\qquad [e_3, e_4] = e_2 +(1+\beta) e_3,\qquad [e_4, e_5] = e_1 + \beta e_2.\end{split}
\end{gather*}
The $\big(\fg^{\langle 1\rangle}\big)_{\ev}$-module structure in~$\big(\fg^{\langle 1\rangle}\big)_{\od}$ is given by the following formulas:
\begin{alignat*}{5}
& [e_1, o_1] = o_1+o_2 +\beta^2 o_4, \qquad && [e_2, o_1] = o_2, \qquad && [e_3, o_1] = o_3+o_4, \qquad && [e_4, o_1] = 0,& \\
& [e_1, o_2] = o_2, \qquad && [e_2, o_2] = 0, \qquad && [e_3, o_2] = o_4, \qquad && & \\
& [e_4, o_2] = o_1 + (1 + \beta) o_2 + \beta ^2 o_3, && && && & \\
& [e_1, o_3] = o_3+o_4, \qquad && [e_2, o_3] = o_4, \qquad && [e_3, o_3] = 0, \qquad && &\\
& [e_4, o_3] = o_2 + \beta (o_3+ o_4), && && && & \\
& [e_1, o_4] = o_4, \qquad && [e_2, o_4] = 0, \qquad && [e_3, o_4] = 0, \qquad && [e_4, o_4] = o_3+o_4,& \\
& [e_5, o_1] = o_2 + \beta^2 o_4, \qquad && [e_5, o_2] = 0, \qquad && [e_5, o_3] = o_4, \qquad && [e_5, o_4] = 0.&
\end{alignat*}

{\it Grading $c=(0,0,1)$.} The Lie subsuperalgebra $\fg^{\langle 1\rangle}$, see equation~\eqref{subSup}, is as follows:
\begin{gather*}
 \big(\fg^{\langle 1\rangle}\big)_{\ev} = \Span\big(e_1 = \big(x + x^{(4)}\big)\del,\, e_2 = x^{(5)}\del,\\
\hphantom{\big(\fg^{\langle 1\rangle}\big)_{\ev} = \Span\big(}{} e_3 = \big(x^{(3)} + x^{(6)}\big)\del,\, e_4 = \del^2 + x^{(2)}\del,\, e_5 = x^{(7)}\del\big),\\
 \big(\fg^{\langle 1\rangle}\big)_{\od} = \Span\big(o_1 = \big(1+x^{(3)}\big)\del,\, o_2 = x^{(4)}\del,\, o_3 = \big(x^{(2)} + x^{(5)}\big)\del,\, o_4 = x^{(6)}\del\big).
\end{gather*}
The squaring in $\fs\big(\fg^{\langle 1\rangle}\big)$ is given by the following formulas:
\begin{gather*}
 o_1^2 = e_4,\qquad o_2^2 = e_5,\qquad o_3^2 = e_3,\qquad o_4^2= 0.
\end{gather*}
The non-zero commutation relations in~$\big(\fg^{\langle 1\rangle}\big)_{\ev}$ are given by the following formulas:
\begin{gather*}
[e_1, e_4]= e_2,\qquad [e_2, e_4]= e_3,\qquad [e_3, e_4]= e_1,\qquad [e_4, e_5]=e_2.
\end{gather*}
The $\big(\fg^{\langle 1\rangle}\big)_{\ev}$-module structure in~$\big(\fg^{\langle 1\rangle}\big)_{\od}$ is given by the following formulas:
\begin{alignat*}{6}
& [e_1, o_1] = o_1+o_4, \qquad && [e_2, o_1] = o_2, \qquad && [e_3, o_1] = o_3, \qquad && [e_4, o_1] = 0, \qquad && [e_5, o_1] = o_4,& \\
& [e_1, o_2] = o_2, \qquad && [e_2, o_2] = 0, \qquad && [e_3, o_2] = o_4, \qquad && [e_4, o_2] = o_3, \qquad && [e_5, o_2] = 0,& \\
& [e_1, o_3] = o_3, \qquad && [e_2, o_3] = o_4, \qquad && [e_3, o_3] = 0, \qquad && [e_4, o_3] = o_1+o_4, \qquad && [e_5, o_3] = 0, & \\
& [e_1, o_4] = o_4, \qquad && [e_2, o_4] = 0, \qquad && [e_3, o_4] = 0, \qquad && [e_4, o_4] = o_2, \qquad && [e_5, o_4] = 0. &
\end{alignat*}

{\it Grading $c=(1,0,0)$.} The Lie subsuperalgebra $\fg^{\langle 1\rangle}$, see equation~\eqref{subSup}, is as follows:
\begin{gather*}
 \big(\fg^{\langle 1\rangle}\big)_{\ev} = \Span\big(e_{-2} = \del^2 + (1+x)\del,\, e_{-1} = (1+x)\del,\, e_{0} = \big(x+x^{(2)}+x^{(3)}\big)\del,\\
\hphantom{\big(\fg^{\langle 1\rangle}\big)_{\ev} = \Span\big(}{} e_{1} = \big(x^{(2)}+x^{(4)}+x^{(5)}\big)\del,\, e_{2} = \big(x^{(3)}+x^{(5)}+x^{(6)}+x^{(7)}\big)\del\big),\\
 \big(\fg^{\langle 1\rangle}\big)_{\od} = \Span\big(o_1 = \del,\, o_2 = \big(x+x^{(2)}\big)\del,\\
\hphantom{\big(\fg^{\langle 1\rangle}\big)_{\od} = \Span\big(}{} o_3=\big(x^{(2)}+x^{(3)}+x^{(4)}\big)\del,\, o_4 = \big(x^{(3)}+x^{(5)}+x^{(6)}\big)\del\big).
\end{gather*}
The squaring in $\fs\big(\fg^{\langle 1\rangle}\big)$ is given by the following formulas:
\begin{gather*}
 o_1^2 = e_{-2}, \qquad o_2^2 = e_0, \qquad o_3^2 = e_2,\qquad o_4^2 = 0.
\end{gather*}
The non-zero commutation relations in~$\big(\fg^{\langle 1\rangle}\big)_{\ev}\simeq(\fvect(1;\underline{3}))^{\langle 1\rangle}$ are given by the following formulas:
\begin{gather*}
[e_{-1}, e_0]= e_{-1}, \qquad [e_{-1}, e_1]= e_0, \qquad [e_{-1}, e_2]= e_1, \qquad [e_{1}, e_{-2}]= e_{-1}, \qquad [e_0,e_1]= e_1.
\end{gather*}
The $\big(\fg^{\langle 1\rangle}\big)_{\ev}$-module structure in~$\big(\fg^{\langle 1\rangle}\big)_{\od}$ is given by the following formulas:
\begin{alignat*}{5}
& [e_{-1}, o_1] = o_1, \qquad && [e_0, o_1] = o_1 + o_2, \qquad && [e_1, o_1] = o_2 + o_3, \qquad && [e_2, o_1] = o_3 + o_4, &\\
& [e_{-2}, o_1] = o_1, &&&&&& &\\
& [e_{-1}, o_2] = o_1 + o_2, \qquad && [e_0, o_2] = 0, \qquad && [e_1, o_2] = o_3 + o_4, \qquad && [e_2, o_2] = 0, &\\
& [e_{-2}, o_2] = o_2, &&&&&& & \\
& [e_{-1}, o_3] = o_2 + o_3, \qquad && [e_0, o_3] = o_3 + o_4, \qquad && [e_1, o_3] = 0, \qquad && [e_2, o_3] = 0, &\\
& [e_{-2}, o_3] = o_1 + o_3, &&&&&& & \\
& [e_{-1}, o_4] = o_3 + o_4, \qquad && [e_0, o_4] = 0, \qquad && [e_1, o_4] = 0, \qquad && [e_2, o_4] = 0, &\\
& [e_{-2}, o_4] = o_2 + o_4.&&&&&& &
\end{alignat*}
We have
\begin{gather*}
\big[\big(\fg^{\langle 1\rangle}\big)_{\ev}, \big(\fg^{\langle 1\rangle}\big)_{\ev}\big] = \Span(e_{-1},\, e_0,\, e_1),\\
\big[\big(\fg^{\langle 1\rangle}\big)_{\ev}, \big[\big(\fg^{\langle 1\rangle}\big)_{\ev}, \big(\fg^{\langle 1\rangle}\big)_{\ev}\big]\big]
= \big[\big(\fg^{\langle 1\rangle}\big)_{\ev}, \big(\fg^{\langle 1\rangle}\big)_{\ev}\big].
\end{gather*}
The ideal $[\big(\fg^{\langle 1\rangle}\big)_{\ev}, \big(\fg^{\langle 1\rangle}\big)_{\ev}]$ is isomorphic to~${\mathfrak{o}}^{(1)}(3)$. The ${\mathfrak{o}}^{(1)}(3)$-module $\big(\fg^{\langle 1\rangle}\big)_{\od}$ is irreducible, it has no lowest weight vectors and has highest vectors~$o_3$ and~$o_4$, cf.~\cite{Do}.

\subsection{Summary of computer-aided experiments}\label{SSSumComp}
Let $L_k$ and $L^{(k)}$ be given by equation~\eqref{L_k}; in tables below, $k=1,2,\dots$ up to first stable term. Theorem~\ref{T_vect} implies the following fact:
\begin{gather*}
\sdim \fs\big(\fg, D_{u(c)}\big)=\big(2^{n-1}+1|2^{n-1}\big)\qquad \text{and}\\ \dim\big[\big(\fg^{\langle 1\rangle}\big)_{\ev}, \big(\fg^{\langle 1\rangle}\big)_{\ev}\big]=2^{n-1}-1.
\end{gather*}
Recall that $\fv_n$ is defined by~\eqref{fv_n}; the symbol``\textsf{solv}'' below means that the corresponding Lie algebra is solvable, see Corollary~\ref{Cor_Solvable}. Computer-aided experiments show that \begin{gather*}
\fvect^{(1)}(1;\underline{2})\quad \renewcommand{\arraystretch}{1.2}
\begin{tabular}{|@{\,}l@{\,}|@{\,}l@{\,}|@{\,}l@{\,}|@{\,}l@{\,}|} \hline
 $\big(\fg^{\langle 1\rangle}\big)_{\ev}$ & $\dim L_{k}$ & $\dim L^{(k)}$ & parameters \\[1mm]
 \hline
 Heisenberg & $1,0$ & $1,0$ & $(00)$\\
 \hline
 \textsf{solv}, see \eqref{neHei} & $1$ & $1,0$ & $(01)$, $(10)$, $(11)$ \\
 \hline
\end{tabular} \\[1ex]
\fvect^{(1)}(1;\underline{3})\quad \renewcommand{\arraystretch}{1.2}
\begin{tabular}{|@{\,}l@{\,}|@{\,}l@{\,}|@{\,}l@{\,}|@{\,}l@{\,}|} \hline
 $\big(\fg^{\langle 1\rangle}\big)_{\ev}$ & $\dim L_{k}$ & $\dim L^{(k)}$ & parameters \\[1mm]
 \hline
 ${\mathfrak{o}}(3)/\fc\simeq\fv_3$ & $3$ & $3$ & $(1ab)$ \\[1mm]
 \hline
 \textsf{solv} & $3,2,1,0$ & $3,0$ & $(000)$\\[1mm]
 \hline
 \textsf{solv} & $3,2$ & $3,0$ & $(010)$, $(010)$, $(0\alpha0)$,\text{~where $\alpha\neq0$} \\[1mm]
 \hline
 \textsf{solv} & $3$ & $3,0$ & $(001)$, $(011)$, $(0\alpha\beta)$,\text{~where~$\beta\neq0$}\\[1mm]
 \hline
\end{tabular}
\\[1ex]
\fvect^{(1)}(1;\underline{4})\quad{ \renewcommand{\arraystretch}{1.2}
\begin{tabular}{|@{\,}l@{\,}|@{\,}l@{\,}|@{\,}l@{\,}|@{\,}l@{\,}|} \hline
 $\fg_{\ev}^{\langle 1\rangle}$ & $\dim L_k$ & $\dim L^{(k)}$ & parameters \\[1mm]
 \hline
 $\fv_4$ & $7$ &$7$ &
 $(1abc)$\\
 \hline
 \textsf{solv} & $7,6,5,4,3,2,1,0$ & $7,0$ & $(0000)$ \\[1mm]
 \hline
 \textsf{solv} & $7,6,4$ & $7,0$ & $(0100)$ \\[1mm]
 \hline
 \textsf{solv} & $7,6$ & $7,0$ & $(0010)$, $(0110)$ \\[1mm]
 \hline
 \textsf{solv} & $7$ & $7,0$ & $(0001)$, $(0011)$, $(0101)$, $(0111)$ \\[1mm]
 \hline
\end{tabular}}
\\[1ex]
\fvect^{(1)}(1;\underline{5})\quad{ \renewcommand{\arraystretch}{1.2}
\begin{tabular}{|@{\,}l@{\,}|@{\,}l@{\,}|@{\,}l@{\,}|@{\,}l@{\,}|} \hline
 $\fg_{\ev}^{\langle 1\rangle}$ & $\dim L_k$ & $\dim L^{(k)}$ & parameters \\[1mm]
 \hline
 $\fv_5$ & $15$ &$15$ &
 $(1abcd)$\\
 \hline
 \textsf{solv} & $15,14,13,12,11,10,9,8,$ & $15,0$ & $(00000)$ \\
 {} & $7,6,5,4,3,2,1,0$ & {} & {} \\[1mm]
 \hline
 \textsf{solv} & $15,14,13,12,11,10,9,8$ & $15,0$ & $(01000)$ \\[1mm]
 \hline
 \textsf{solv} & $15,14,13,12$ & $15,0$ & $(00100)$, $(01100)$ \\[1mm]
 \hline
 \textsf{solv} & $15,14$ & $15,0$ & $(00010)$, $(00110)$, $(01010)$, $(01110)$ \\[1mm]
 \hline
 \textsf{solv} & $15$ & $15,0$ & $(00001)$, $(00011)$, $(00101)$, $(00111)$, \\
 {} & {} & {} &$(01001)$, $(01011)$, $(01101)$, $(01111)$ \\[1mm]
 \hline
\end{tabular}}
\end{gather*}

Let $\fii$ be the last term of the sequence~$L_k$, see \eqref{L_k}, for the even part~$\fg_{\ev}$ of the corresponding Lie superalgebra~$\fg$. Set
\begin{gather*}
 c(i) := (0\ldots010\ldots0)\qquad \text{with a 1 in the $(i+1)$st slot}.
\end{gather*}
\begin{Conjecture}Grading $c(k)$ with the corresponding generating function $u=x^{(2^k)}$, where $k=1,\ldots,n-1$, yields the Lie superalgebras such that $\dim\fii=2^{n-1} - 2^{n-1-k}$.
\end{Conjecture}

The following two conjectures concern $\fvect^{(1)}(1;\underline{n})$, where $n\geq 3$.

\begin{Conjecture}Let $u = \sum\limits_{1\leq i\leq n-1} c_i x^{(2^i)}$, i.e., $c_0=0$. If all $c_i$ are equal to $0$, then let $k$ be $0$; otherwise, let $k$ be the maximal number such that $c_k\neq0$.
Then the grading given by $D_u$ yields the Lie superalgebras such that $\dim\fii=2^{n-1} - 2^{n-1-k}$.
\end{Conjecture}

\subsection[The derivations of~$\fvect^{(1)}(1;\underline{n})$]{The derivations of~$\boldsymbol{\fvect^{(1)}(1;\underline{n})}$}
The result of this subsection is probably known, but we'd like to draw attention to Sierpi\'nski sieves here. Observe that solutions of the linear equations~\eqref{cond2lin}, i.e., derivations of~$\fvect^{(1)}(1;\underline{n})$, form the Sierpi\'nski sieve of order~$\underline{n}$ under the main diagonal, i.e., for
$i\geq j$, we have
\begin{gather*}
 c_{i,j} = \binom{i}{j-1} c_{i-j+1,1},\qquad \text{which corresponds to the derivation~$\ad_{x^{(i-j+1)}\del}$},
 \end{gather*}
and for $i < j$ we have nonzero diagonals with parameters~$c_{1, 2^k+1}$, where $k=0,\dots,n-1$, which correspond to the derivations~$\ad_{\del^{2^k}}$, e.g.,
\begin{gather*}
 c_{1,j} = 0\qquad\text{if $j\neq 2^k+1$, where $k$ is a non-negative integer},\\
 c_{i,j} = c_{i-1,j-1} \qquad \text{for $i>1$}.
\end{gather*}

Namely, for $\underline{n}=3$, we have the following solution of linear equations~\eqref{cond2lin}, and its schematic form (here $\ast$ represents any nonzero entry and the empty entries represent zeros):
\begin{gather*}
U_{{\rm lin}} =
\begin{pmatrix}
 c_{1,1} & c_{1,2} & c_{1,3} & 0 & c_{1,5} & 0 & 0 \\
 c_{2,1} & 0 & c_{1,2} & c_{1,3} & 0 & c_{1,5} & 0 \\
 c_{3,1} & c_{2,1} & c_{1,1} & c_{1,2} & c_{1,3} & 0 & c_{1,5} \\
 c_{4,1} & 0 & 0 & 0 & c_{1,2} & c_{1,3} & 0 \\
 c_{5,1} & c_{4,1} & 0 & 0 & c_{1,1} & c_{1,2} & c_{1,3} \\
 c_{6,1} & 0 & c_{4,1} & 0 & c_{2,1} & 0 & c_{1,2} \\
 c_{7,1} & c_{6,1} & c_{5,1} & c_{4,1} & c_{3,1} & c_{2,1} & c_{1,1} \\
\end{pmatrix},
\qquad
\begin{pmatrix}
 * & * & * & & * & & \\
 * & & * & * & & * & \\
 * & * & * & * & * & & * \\
 * & & & & * & * & \\
 * & * & & & * & * & * \\
 * & & * & & * & & * \\
 * & * & * & * & * & * & * \\
\end{pmatrix}.
\end{gather*}

For $\underline{n}=4$, we have the following schematic form of the solution of linear equations~\eqref{cond2lin}:
\begin{gather*}
\left(
\begin{array}{ccccccccccccccc}
 * & * & * & & * & & & & * & & & & & & \\
 * & & * & * & & * & & & & * & & & & & \\
 * & * & * & * & * & & * & & & & * & & & & \\
 * & & & & * & * & & * & & & & * & & & \\
 * & * & & & * & * & * & & * & & & & * & & \\
 * & & * & & * & & * & * & & * & & & & * & \\
 * & * & * & * & * & * & * & * & * & & * & & & & * \\
 * & & & & & & & & * & * & & * & & & \\
 * & * & & & & & & & * & * & * & & * & & \\
 * & & * & & & & & & * & & * & * & & * & \\
 * & * & * & * & & & & & * & * & * & * & * & & * \\
 * & & & & * & & & & * & & & & * & * & \\
 * & * & & & * & * & & & * & * & & & * & * & * \\
 * & & * & & * & & * & & * & & * & & * & & * \\
 * & * & * & * & * & * & * & * & * & * & * & * & * & * & * \\
\end{array}
\right)\end{gather*}

\appendix

\section{Necessary proofs from~\cite{BGLL2}}\label{appendixA}

\begin{Statement}[{\cite[Statement 3.8.1.a]{BGLL2}}] For $n=3$ or $n\geq 5$, the algebra $\fder\, {\mathfrak{o}}^{(1)}_I(n)$ can be identified with ${\mathfrak{o}}_I(n)/\fc$ in the sense that for any $D\in\fder\, {\mathfrak{o}}^{(1)}_I(n)$, there is $A_D\in{\mathfrak{o}}_I(n)$ such that~$D$ coincides with the restriction of $\ad_{A_D}$ to ${\mathfrak{o}}^{(1)}_I(n)$; for a given $D$, the element $A_D$ is uniquely defined up to adding a~scalar matrix.
\end{Statement}

\begin{proof} In this proof, $i$, $j$, $k$, $l$, $m$ are always indices from $1$ to $n$.

The algebra ${\mathfrak{o}}^{(1)}_I(n)$ consists of zero-diagonal symmetric $n\times n$ matrices, which means that the elements $\bo{i}{j} := E^{i,j}+E^{j,i}$, where $\{i,j\}$ are all two-element subsets of $\{1,\dots,n\}$, form a~basis of ${\mathfrak{o}}^{(1)}_I(n)$. Their commutation relations are (we assume that $i\neq j$ and $k\neq l$):
\begin{gather*}
{}[\bo{i}{j}, \bo{k}{l}] = \begin{cases}0&\text{if~} \{k,l\} = \{i,j\} \text{~or~} \{k,l\} \cap \{i,j\} = \varnothing,\\
 \bo{j}{k}&\text{for~} l=i\text{~and~}k\neq i, j.
 \end{cases}
\end{gather*}

Alternatively, we can say that for an arbitrary matrix $M\in {\mathfrak{o}}^{(1)}_I(n)$ and $i\neq j$,
\begin{gather*}
\big[M, \bo{i}{j}\big]_{kl} = 0 \qquad \text{if} \ \{k,l\} = \{i,j\} \text{~or~} \{k,l\} \cap \{i,j\} = \varnothing,\\
\big[M, \bo{i}{j}\big]_{ik} = \big[M, \bo{i}{j}\big]_{ki} = M_{jk} \qquad \text{for~} k\neq i,j,\\
\big[M, \bo{i}{j}\big]_{kk} = 0 \qquad \text{for an arbitrary~}k.
\end{gather*}

Let $D$ be a derivation of ${\mathfrak{o}}^{(1)}_I(n)$. Let us prove that for arbitrary three pairwise distinct indices~$i$, $j$, $k$,
\begin{gather}\label{oI-eq3}
\big(D\bo{i}{j}\big)_{ij} + \big(D\bo{i}{k}\big)_{ik} + \big(D\bo{j}{k}\big)_{jk} = 0,
\end{gather}
and that for arbitrary four pairwise distinct indices $i$, $j$, $k$, $l$,
\begin{gather} \label{oI-eq1}
 \big(D\bo{i}{j}\big)_{kl} = 0,\\
 \label{oI-eq2}
 \big(D\bo{i}{k}\big)_{il} = \big(D\bo{j}{k}\big)_{jl} = \big(D\bo{i}{l}\big)_{ik}.
\end{gather}

Since $\bo{i}{k} = \big[\bo{i}{j}, \bo{j}{k}\big]$, we have
\begin{gather*}
\big(D\bo{i}{k}\big)_{ik} = \big(D[\bo{i}{j}, \bo{j}{k}]\big)_{ik} = \big[D\bo{i}{j}, \bo{j}{k}\big]_{ik} + \big[\bo{i}{j}, D\bo{j}{k}\big]_{ik}\\
\hphantom{\big(D\bo{i}{k}\big)_{ik}}{} = \big(D\bo{i}{j}\big)_{ij} + \big(D\bo{k}{j}\big)_{jk},
\end{gather*}
which proves (\ref{oI-eq3}).

Note that (\ref{oI-eq1}) and (\ref{oI-eq2}) are vacuously true for $n=3$, since in this case, there exist no four pairwise distinct indices. In case of $n\geq 5$, let $m$ be an index different from all of $i$, $j$, $k$, $l$. Then $\big[\bo{i}{j}, \bo{l}{m}\big] = 0$, which means that
\begin{gather*}
0 = \big(D\big[\bo{i}{j}, \bo{l}{m}\big]\big)_{km} = \big[D\bo{i}{j}, \bo{l}{m}\big]_{km} + \big[\bo{i}{j}, D\bo{l}{m}\big]_{km} = \big(D\bo{i}{j}\big)_{kl} + 0,
\end{gather*}
which proves (\ref{oI-eq1}). Since $\bo{i}{k} = \big[\bo{i}{j}, \bo{j}{k}\big]$, we have
\begin{gather*}
\big(D\bo{i}{k}\big)_{il} = \big(D\big[\bo{i}{j}, \bo{j}{k}\big]\big)_{il} = \big[D\bo{i}{j}, \bo{j}{k}\big]_{il} + \big[\bo{i}{j}, D\bo{j}{k}\big]_{il} = 0 + \big(D\bo{j}{k}\big)_{jl},
\end{gather*}
which proves the first equality in (\ref{oI-eq2}). And since $\bo{i}{k} = \big[\bo{i}{l}, \bo{k}{l}\big]$, we have{\samepage
\begin{gather*}
\big(D\bo{i}{k}\big)_{il} = \big(D\big[\bo{i}{l}, \bo{k}{l}\big]\big)_{il} = \big[D\bo{i}{l}, \bo{k}{l}\big]_{il} + \big[\bo{i}{l}, D\bo{k}{l}\big]_{il} = \big(D\bo{i}{l}\big)_{ik} + 0,
\end{gather*}
which proves the second equality in (\ref{oI-eq2}).}

Now consider matrix $A_D$ whose entries are as follows:
\begin{gather*}
 (A_D)_{11} = 0,\qquad
 (A_D)_{ii} = \big(D\bo{1}{i}\big)_{1i} \qquad \text{for} \ i\neq 1,\\
 (A_D)_{ij} = \big(D\bo{k}{i}\big)_{kj} \qquad \text{for} \ i\neq j, \text{~and some~} k\neq i,j.
\end{gather*}
Note that (\ref{oI-eq2}) shows that $(A_D)_{ij}$ does not depend on the choice of $k$ and that the matrix is symmetric, i.e., $(A_D)_{ij} = (A_D)_{ji}$.

Then it follows from (\ref{oI-eq3}), (\ref{oI-eq1}) and (\ref{oI-eq2}) that $D\bo{i}{j} = \big[A_D, \bo{i}{j}\big]$ for an arbitrary $\bo{i}{j}$. Since the $\bo{i}{j}$ form a basis of ${\mathfrak{o}}^{(1)}_I(n)$, it means that $D$ coincides with the restriction of $\ad_{A_D}$ to ${\mathfrak{o}}^{(1)}_I(n)$. On the other hand, it is easy to check that for any matrix $A\in{\mathfrak{o}}_I(n)$, the restriction of $\ad_{A}$ to ${\mathfrak{o}}^{(1)}_I(n)$ is a derivation of ${\mathfrak{o}}^{(1)}_I(n)$, and that two matrices $A, A'\in{\mathfrak{o}}_I(n)$ determine the same derivation of ${\mathfrak{o}}^{(1)}_I(n)$ if and only if $A-A' = c1_n$ for some $c\in {\mathbb K}$. This proves our claim.
\end{proof}

\subsection*{Acknowledgements} We are thankful to D.~Leites, who raised the problem, I.~Shche\-poch\-kina, and S.~Skryabin for help, and to P.~Grozman, whose code \textit{SuperLie}, see \cite{Gr}, we used in our computer experiments. Thanks are due to S.~Bouarroudj, and V.~Grandjean for discussions and useful comments. The first author thanks the Organising committee of the symposium ``Groningen Deformation Day'' (October~7, 2016, Groningen, The Netherlands), where the results of this note were delivered, for hospitality and financial support; his research was partly supported by WCMCS post-doctoral fellowship and the grant AD 065 NYUAD during his visits of NYUAD. For the possibility to perform the difficult computations of this research we are grateful to M.~Al~Barwani, Director of the High Performance Computing resources at New York University Abu Dhabi.

\pdfbookmark[1]{References}{ref}
\LastPageEnding

\end{document}